\documentclass[12pt]{article}

\def\qed{\quad \vrule height7.5pt width4.17pt depth0pt}
 
\newtheorem{theorem}{Theorem}
\newtheorem{lemma}{Lemma}[section]
\newtheorem{prop}{Proposition}[section]

\newtheorem{definition}{Definition}[section]

\begin{document}

\title{Extinction for two parabolic stochastic PDE's on the lattice}
\author{C. Mueller$^1$ and E. Perkins$^2$\\\\
  Dept. of Mathematics\\
  University of Rochester\\
  Rochester, NY  14627\\
  USA\\
  E-mail:  cmlr@troi.cc.rochester.edu\\\\
  Dept. of Mathematics\\
  University of British Columbia\\
  Vancouver, BC  V6T 1Z2\\
  Canada\\
  E-mail:  perkins@heron.math.ubc.ca}
\date{}
\maketitle
\footnotetext[1]{Research supported in part by an NSA grant and
an NSERC Collaborative Grant}
\footnotetext[2]{Research supported in part by an NSERC Research
Grant and an NSERC Collaborative Grant


AMS 1991 {\em subject classifications}
Primary, 60H15; Secondary, 35R60, 35L05.}

\newpage

\begin{abstract}
It is well known that, starting with finite mass, the super-Brownian
motion dies out in finite time.  The goal of this article is to show
that with some additional work, one can show finite time die-out
for two types of systems of stochastic differential equations 
on the lattice $\mathbf{Z}^d$.  

For our first system, let $1/2\le\gamma<1$, and consider non-negative 
solutions of  
\begin{eqnarray*}
du(t,x)&=&\Delta u(t,x)dt+u^\gamma(t,x) dB_x(t),\ x\in \mathbf{Z}^d	\\
u(0,x)&=&u_0(x)\ge 0.	
\end{eqnarray*}
Here $\Delta$ is the discrete Laplacian and $\{B_x:x\in \mathbf{Z}^d\}$ 
is a system of independent Brownian motions.  
We assume that $u_0$ has finite support.
When $\gamma=1/2$, the measure which puts mass $u(t,x)$ at $x$ is a 
super-random walk and it is well-known that the process becomes extinct
in finite time a.s.  Finite-time extinction is known to be a.s. false if 
$\gamma=1$.
For $1/2<\gamma<1$, we show finite-time
die-out by breaking up the solution into pieces, and showing that each
piece dies in finite time.  Unlike the superprocess case, these
pieces will not in general evolve independently.

Our second example involves the mutually catalytic branching 
system of stochastic differential equations on $\mathbf{Z^d}$, which was
first studied in Dawson and Perkins \cite{dp98}. 
\begin{eqnarray*}
dU_t(x)&=&\Delta U_t(x)dt+\sqrt{U_t(x)V_t(x)} dB_{1,x}(t)	\\
dV_t(x)&=&\Delta V_t(x)dt+\sqrt{U_t(x)V_t(x)} dB_{2,x}(t)	\\
U_0(x)&\ge& 0							\\
V_0(x)&\ge& 0.	
\end{eqnarray*}
By using a somewhat different argument, we show that, depending on the
initial conditions, finite time extinction of one type may occur with 
probability 0, or with probability arbitrarily close to 1.  
\end{abstract}

\newpage

\section{Introduction}
  \setcounter{equation}{0}

Recently, the Dawson-Watanabe process, or super-Brownian motion, has
attracted great interest, and many fascinating properties have come to
light.  See Dawson \cite{daw93} for a survey.  These results often rely
on the multiplicative property of the process.  This allows one to study
the process as an infinite system of noninteracting 
particles, each with infinitesimal mass.  However, it is often much more
difficult to prove similar results for systems with interactions.  

In this article, we concentrate on the finite time extinction property.
Let $Z_t$ be the total mass of the Dawson-Watanabe process, and assume
that the initial mass $Z_0<\infty$.  As is well known, $Z_t$ satisfies
the Feller equation
\begin{equation}
\label{fell}
dZ=\sqrt{Z}dB
\end{equation}
and with probability 1, $Z_t$ reaches 0 in finite time.  See Theorem 4.3.6
of \cite{kni81} for the exact extinction probabilities.
Our goal is to study finite time extinction
for 2 types of systems of stochastic differential equations (SDE), related
to super-random walks,  
on the lattice 
$\mathbf{Z}^d$.  

First, we consider non-negative solutions $u(t,x)$, $ t\ge 0,x\in\mathbf{Z}^d$
to the following system of stochastic differential equations on the 
lattice $\mathbf{Z}^d$, for $1/2\le\gamma<1$.
\begin{eqnarray}
\label{system}
du(t,x)&=&\Delta u(t,x)dt+u^\gamma(t,x) dB_x(t),\ x\in {\bf Z}^d	\\
u(0,x)&=&u_0(x)\ge 0.	\nonumber
\end{eqnarray}
Here and throughout the paper, $\Delta$ is the discrete Laplacian on
$\mathbf{Z}^d$.  In other words, if
$\mathcal{N}(x)$ is the set of $2d$ nearest neighbors of
$x\in\mathbf{Z}^d$, and if $f(x)$ is a function on $\mathbf{Z}^d$, then
$(\Delta f)(x)=\sum_{y\in\mathcal{N}(x)}f(y)-2df(x)$.
Also, $\{B_x(t)\}_{x\in\mathbf{Z}^d}$ is a collection of
independent $(\mathcal{F}_t)$-Brownian motions on some filtered
probability space $(\Omega,\mathcal{F},\mathcal{F}_t,P)$ satisfying the
usual right-continuity and completion hypotheses, as will all our
filtered probability spaces in this work.    We assume that
$u_0(x)$ equals 0 except at a finite number of points, $\mathbf{F}$, 
in $\mathbf{Z}^d$.  
Pathwise existence and uniqueness holds for solutions of (\ref{system}) by 
the well-known method of Yamada and Watanabe which we recall below (Lemma
\ref{unique}).  

If $\gamma =1$, solutions to (\ref{system}) can be represented in
terms of the Feynman-Kac formula.  Let $\xi(t)$ be a continuous time
random walk on $\mathbf{Z}^d$ with infinitesimal generator $\Delta$ and semigroup
$P_t$, which
is independent of the Brownian motions $B_x$.  If $\mathbf{E}_x$ denotes
the expectation with respect to $\xi$, for $\xi(0)=x$, then we have
\[
u(t,x)=\mathbf{E}_x\left(u_0(\xi(t))
	\exp\left[\int_0^tdB_{\xi(t-s)}(s)-t/2\right]\right).
\]
Since $\exp[\cdot]$ is always strictly positive, and since for each
$t>0$ there is a positive probability that $\xi(t)$ lies in the support of
$u_0$, it follows that $u(t,x)>0$ for all $t>0$, $x\in\mathbf{Z}^d$.
G\"artner and Molchanov \cite{gm90} have found many fascinating
properties of solutions for the case $\gamma=1$.  We also mention in
passing that a class of processes called ``linear systems'' has been studied
in the particle systems literature.  Such systems are formally similar to
solutions of (\ref{system}) with $\gamma=1$, and Liggett, \cite{lig85}
gives some theorems about the asymptotic die-out of mass as $t\to\infty$.  

Next, for $\gamma=1/2$, the measure $u(t,x)d\mu(x)$, (where $\mu$ is the
counting measure), is the super-Brownian motion with underlying spatial
motion $\xi(t)$.  Its total mass satisfies (\ref{fell}) and 
therefore, $u(t,x)=0$ for all $x$ and large enough $t$

In light of the above two results it is natural to consider the question
of finite time extinction for $1/2<\gamma<1$.  In this case one can view 
solutions to (\ref{system}) as interactive super-random walks in which
there is a density dependent branching rate of $u(t,x)^{\gamma-1/2}$ at
$(t,x)$.  Clearly, for some Brownian motion $B(t)$, the total mass $Z(t)$
satisfies 
\begin{equation}
\label{feller}
dZ(t)=\left(\sum_{x\in\mathbf{Z}^d}u^{2\gamma}(t,x)\right)^{1/2}dB(t).
\end{equation}
Suppose that $H(t)\ge cZ^{\gamma}(t)$, where $H(t)$ is nonanticipating. 
It is known that with probability 1, solutions to $dZ=HdB$ die out in
finite time.  See, for example, Lemma 3.4 of \cite{mp92}.  Unfortunately,
if $u(t,x)$ is very thinly spread, 
$[\sum_{x\in\mathbf{Z}^d}u^{2\gamma}(t,x)]^{1/2}$
may be much smaller than $Z^{\gamma}(t)$.  Thus, the coefficient of
$dB(t)$ which appears in (\ref{feller}) may be much smaller than 
$Z^{\gamma}(t)$.  This is the main difficulty in proving Theorem 1.  

\begin{theorem}
\label{t1}
Suppose that $1/2\le\gamma<1$, $d\ge 1$, that $u(t,x)$ satisfies
(\ref{system}), and that $u_0(x)$ is equal to 0 except on a finite set 
$\mathbf{F}$.  Then, with probability 1, $u(t,x)$ dies out in finite time.
That is, there exists an almost surely finite random time
$\tau=\tau(\omega)$ such that $u(t,x)=0$ for all $t\ge\tau$ and 
$x\in\mathbf{Z}^d$.   
\end{theorem}

The strategy of our proof is to show that $u(t,x)$ is not thinly spread,
and therefore $Z(t)$ satisfies an equation like $dZ=HdB$, where 
\begin{equation}
\label{clumping}
H\ge cZ^\gamma(t) 
\end{equation}
for some random number $c$.  It is known that solutions to such equations
die out in finite time.  Actually, $u(t,x)$ shows a high degree of
clumping as $x$ varies.  Here, we were guided by known results for
superprocesses.  Without the clumping, we would not be able to show an
inequality such as (\ref{clumping}).

Finite-time extinction is often useful for establishing the compact
support property of solution to continuous parameter parabolic stochastic
PDE's.  The compact support property states that if the initial data has
compact support in $\mathbf{R}$ then the same is true of the solution at
any positive time.
Often, the compact support property of solutions is proved by showing that
finite-time die-out occurs for the parts of the solution corresponding to large
values of $x$.  For example, this is done in \cite{mp92} and \cite{dp91}.
In fact, Theorem 3.10 of \cite{mp92} proves the continuous analogue of 
Theorem \ref{t1} but we were unable to extend that approach to our lattice
systems.  At a crucial step in the proof in \cite{mp92}, we used Jensen's
inequality.  To prove Theorem \ref{t1} we again use Jensen's inequality,
but we also need to know that the mass of $u(t,x)$ tends to
cluster at a small number of sites.   

Next, we introduce a system of SDE's introduced in \cite{dp98}.  
Let $M_F(\mathbf{Z}^d)$ be the 
space of finite measures on $\mathbf{Z}^d$ with the topology of weak
convergence.  Consider
\begin{eqnarray}
\label{system2}
dU_t(x)&=&\Delta U_t(x)dt+\sqrt{U_t(x)V_t(x)} dB_{1,x}(t), x\in \mathbf{Z}^d
\nonumber	\\ dV_t(x)&=&\Delta V_t(x)dt+\sqrt{U_t(x)V_t(x)} dB_{2,x}(t), x\in
\mathbf{Z}^d
 \\
U_0, V_0&\in& M_F(\mathbf{Z}^d).				\nonumber	
\end{eqnarray}
Here, $\{B_{i,x}(t)\}_{x\in\mathbf{Z}^d; i=1,2}$ is a collection of independent
$\mathcal{F}_t$-Brownian motions ($\mathcal{F}_t$ are as above)
and $\Delta$ is the discrete Laplacian on $\mathbf{Z}^d$.

Such a pair of processes arise as the large population limit of two
interacting branching populations in which the branching rate of each
type at $x\in\mathbf{Z}^d$ is proportional to the amount of the other 
type at $x$.  As each type ``catalyzes" the reproduction of the other 
type, it is called the mutually catalytic branching process.  One reason
for interest in this system is that it was an extremely simple example of 
interactive branching for which uniqueness in law was not known.
\cite{myt98} and \cite{dp98} proved
weak existence and uniqueness of solutions to (\ref{system2}) by means of 
a self-duality argument proposed by Mytnik.  Uniqueness in law for general 
systems involving interactive branching rates remains unresolved even
for quite smooth rates.  The original reason for interest in (\ref{system2})
was the qualitative behaviour of its continuum analogues in 2 or more
dimensions (the one dimensional case is treated in \cite{dp98}).  The
singularity of 
super-Brownian motion (for $d\ge 2$) and the fact that each type solves
the heat equation in the absence of the other, suggests that the two types
separate and have densities away from their ``interface".  In
\cite{defmpx99} this description is made precise, at least for $d=2$.  

The components $U_{\cdot},V_{\cdot}$ are continuous $M_F(\bf{Z}^d)$-valued
processes a.s., so let $P_{U_0,V_0}$ denote the law of the solution on 
$C([0,\infty),M_F(\mathbf{Z}^d)^2)$.  We set $\langle U,
\phi\rangle=\sum_{x\in\mathbf{Z}^d}
\phi(x)U(x)$ for bounded $\phi$ and $U\in M_F(\mathbf{Z}^d)$.  
The long time behavior of solutions to (\ref{system2}) was studied in 
\cite{dp98}.  
It is easy to
see that $(\langle U_t,1\rangle,\langle V_t,1\rangle)$ 
is a conformal martingale in
the first quadrant and hence converges a.s. as $t\to\infty$ to 
$(\langle U_\infty,1\rangle,\langle V_\infty,1\rangle)$, 
say.  \cite{dp98} showed that 
in the recurrent case ($d\le 2$), $\langle U_\infty,1\rangle \langle 
V_\infty,1\rangle =0$ a.s.,
while in the transient case ($d\ge 3$), $P_{U_0,V_0}
(\langle U_\infty,1\rangle\langle 
V_\infty,1\rangle >0)>0$.  The self-duality then allowed one to use these 
results to study the long time behaviour from infinite initial conditions.  
For finite initial conditions
the above results lead one to ask:\hfil\break\medskip
\noindent
1. Is there finite-time extinction of one type if $d\le 2$? \hfil\break
\medskip\noindent
2. How large is $P_{U_0,V_0}(\langle U_\infty,1\rangle\langle 
V_\infty,1\rangle >0)$ for
$d\ge 3$?\hfil\break
\medskip\noindent
The next two results show that, depending on the 
initial conditions and independent of the dimension, finite
time extinction can occur with probability zero or probability very close
to one.  In
particular, this shows that in the transient case, $P_{U_0,V_0}(\langle
V_\infty ,1\rangle\langle V_\infty ,1\rangle >0)$ may be arbitrarily small,
depending on $(U_0,V_0)$.

Our first theorem about this system establishes conditions under which 
finite time die-out does not occur.  As usual 
$U_0P_t(x)=\sum_y U_0(y)p_t(y,x)$,
where $p_t(x,y)=P(\xi_t=x|\xi_0=y)$.  If $x=(x_1,\dots,x_d)\in \bf{Z}^d$,
then
$|x|=\sum_{i=1}^d|x_i|$.
\begin{theorem}
\label{3.4}
Assume that for $t$ large enough (say $t>t_0$),
\begin{equation}
\label{3.16}
\liminf_{|x|\to\infty}\frac{U_0P_t(x)}{V_0P_t(x)}=\liminf
_{|x|\to\infty}\frac{V_0P_t(x)}{U_0P_t(x)}=0.
\end{equation}
Then $P_{U_0,V_0}(\langle U_t,1\rangle \langle V_t,1\rangle >0 
	\quad\forall t>0)=1$.
\end{theorem}

In Proposition \ref{3.5} we give a large class of initial conditions
which satisfy (\ref{3.16}).

Our second
result about this case says that under certain conditions, finite time
extinction can occur, at least for one of the species.  

\begin{theorem}
\label{3.7}
Assume $\lambda_1\ge\lambda_2>0$, $\lambda_0>4\lambda_1-3\lambda_2$, and
for some
$0<c_1\leq c_2$,
\[
c_1e^{-\lambda_1 |x|}\leq V_0(x)\leq c_2e^{-\lambda_2|x|}\hbox{  for all
}x.
\]
For any $\varepsilon >0$ and $t_1>0$ $\exists \eta >0$ such that if
$U_0(x)\leq \eta e^{-\lambda_0|x|}$ then
\[
P_{U_0,V_0}\left( U_t\equiv 0\quad\forall t\geq t_1\right)\geq 1-\varepsilon.
\]
\end{theorem}

Note the above conditions are satisfied in particular if
$\lambda_0>\lambda_1=\lambda_2$.  

We will see that to achieve smaller
values of
$\varepsilon$, we must take smaller values of $\eta$.  
Both of the above
results will be stated and proved for more general generators than
$\Delta$ in Section 4.  

If $(U,V)$ solves (\ref{system2}), 
then so does $(cU,cV)$ (with a modified initial 
condition).  
We can take $c=K/\eta$ in the above result to see that if 
$\lambda_i$ are as above, $M\ge 1$, 
and $U_0(x)\le Ke^{-\lambda_0 |x|}$, then for any
$\varepsilon,t_1>0$ there is a $c_1$ sufficiently large so that 
the conclusion of Theorem \ref{3.7} holds whenever 
$$c_1e^{-\lambda_1|x|}\le V_0(x)\le c_1Me^{-\lambda_2|x|}.$$

Here is the plan of our article.  We will first deal with (\ref{system}). 
Section 2 contains some lemmas, including uniqueness for 
(\ref{system}).  In Section 3 we prove Theorem \ref{t1}.  In Section 4
we turn to (\ref{system2}) and prove more general versions of 
Theorems \ref{3.4} and \ref{3.7}.

\section{Some Lemmas}
  \setcounter{equation}{0}
In this section, we prove some preliminary facts and lemmas.
Our first lemma gives uniqueness for (\ref{system}).  
This result follows from the method of Yamada and Watanabe (see Theorem
V.40.1 of
Rogers and Williams \cite{rw87}).  Almost the same proof is given below, 
but we
include it for completeness.

\begin{lemma}
\label{unique}
Suppose that $\sum_{x\in\mathbf{Z}^d}u_0(x)<\infty$.  Then pathwise
uniqueness holds for (\ref{system}).  
\end{lemma}

\noindent
\textbf{Proof.  }
Suppose that $u(t,x)$, $v(t,x)$ are 2 solutions of (\ref{system}).  
An easy application of Fatou's Lemma shows that 
\begin{equation}
\label{rw.0}
E(\sum_xu(t,x)+v(t,x))\le \sum_xu(0,x)+v(0,x)<\infty,
\end{equation}
and in fact with a bit more work one can show equality holds in the above.
Note that $||x|^\gamma-|y|^\gamma|^2\le |x-y|^{2\gamma}$, and since 
$\gamma>1/2$, $\int_{0+}u^{-2\gamma}du=\infty$.  Therefore, proceeding  
as in \cite{rw87}, section V.40, we find that for each $x\in\mathbf{Z}^d$,
\begin{eqnarray}
\label{rw.1}
\lefteqn{\left|u(t,x)-v(t,x)\right|}	\\
	&=&\int_0^t\mbox{sgn}(u(s,x)-v(s,x))
	\left[u^\gamma(s,x)-v^\gamma(s,x)\right]dB_x(s) \nonumber\\
  &&+\int_0^t\mbox{sgn}(u(s,x)-v(s,x))
	\left[\Delta u(s,x)-\Delta v(s,x)\right]ds.  \nonumber
\end{eqnarray}
Note that by the definition of $\Delta$, 
\begin{equation}
\label{rw.2}
\sum_{x\in\mathbf{Z}^d}\left|\Delta u(s,x)-\Delta v(s,x)\right|  
	\le\sum_{x\in\mathbf{Z}^d}(4d)\left|u(s,x)-v(s,x)\right| 
\end{equation}
Taking expectations and summing over $x\in\mathbf{Z}^d$ in (\ref{rw.1}),  
and using (\ref{rw.2}), we find that
\[
\sum_{x\in\mathbf{Z}^d}E\left|u(t,x)-v(t,x)\right|
	\le \int_0^t\sum_{x\in\mathbf{Z}^d}(4d)E\left|u(s,x)-v(s,x)\right|ds.
\]
Thus, (\ref{rw.0}) and 
Gronwall's lemma imply that 
\[
\sum_{x\in\mathbf{Z}^d}E\left|u(t,x)-v(t,x)\right|=0
\]
and Lemma \ref{unique} follows.  \qed

Standard arguments show weak existence of solutions to (\ref{system})
(e.g., as in section 2 of \cite{dp98}), and just as for finite
dimensional SDE (see V.17.1 of \cite{rw87}), this and the above result imply 
pathwise existence for solutions of (\ref{system}) and the uniqueness of
its law on $C([0,\infty),M_F(\mathbf{Z}^d))$.

Let $Y$ be a Poisson random variable with parameter $\lambda$, and
suppose that $H$ is a non-negative integer.  Then, 
by Stirling's formula,
\begin{eqnarray}
\label{poisson}
P(Y\ge H)&=&\sum_{k=0}^\infty\frac{\lambda^{(k+H)}}{(k+H)!}e^{-\lambda}
\nonumber
\\ &\le&\frac{\lambda^H}{H!}\sum_{k=0}^\infty
\frac{\lambda^k}{k!} e^{-\lambda}\nonumber	\\
	&\le&\frac{C\lambda^H}{\sqrt{H}H^He^{-H}}.	
\end{eqnarray}

The next lemma follows from an easy application of It\^o's lemma.  See Ch. 5 of
\cite{wal86} for the result in the more delicate continuum setting.

\begin{lemma}
\label{l.int-eq}
(\ref{system}) is equivalent to the following system of integral equations.
\begin{eqnarray}
\label{int-eq}
u(t,x)&=&\sum_{y\in\mathbf{Z}^d}G(t,x-y)u_0(y) \\
	&&+\int_0^t\sum_{y\in\mathbf{Z}^d}G(t-s,x-y)
      \cdot u^\gamma(s,y)dB_y(s),\ x\in \mathbf{Z}^d,  \nonumber 
\end{eqnarray}
where $G(t,x)$ is the fundamental solution of the discrete heat equation 
on  $\mathbf{Z}^d$.  
\end{lemma}

A standard consequence of pathwise uniqueness is
\begin{lemma}
\label{mark}
(\ref{system}) has the strong Markov property with respect to the 
$\sigma$-fields $\mathcal{F}_t$.
\end{lemma}

If $x\in\mathbf{Z}^d$, let $x_i$ denote the $i^{\mbox{th}}$ component of
the vector $x$.  Fix an integer $N>0$, and let 
$D_N=\{x\in\mathbf{Z}^d:  |x_i|\le N \mbox{ for all $i=1,\dots,d$}\}$.
Note that $D_N$ has at most $(3N)^d$ sites.
Let $\partial D_N$ be the boundary of $D_N$.  In other words, let 
$\partial D_N$ be the set of points in $\mathbf{Z}^d\setminus D_N$ which
are nearest neighbors of some point in $D_{N}$.  Recall that
$x,y\in\mathbf{Z}^d$ are nearest neighbors if the Euclidean distance
between them is 1.  For future use, we denote 
${\overline D}_N=D_N\cup\partial D_N$.  We reserve the notation $u(t,x)$ 
for the unique solution of (\ref{system}). 
\begin{lemma}
\label{decomp1}
Let $v_0(x)$ be a non-negative function supported on $D_N$.
Suppose that for $t>0$ and $x\in\mathbf{Z}^d$, $v(t,x)$ satisfies 
\begin{eqnarray}
\label{cutoff}
dv(t,x)&=&\Delta v(t,x)dt+v^\gamma(t,x)dB_x(t)\quad\mbox{if $x\in D_N$}  \\
v(t,x)&=&0\quad\mbox{if $x\in D_N^c$}  \nonumber\\
v(0,x)&=&v_0(x)\quad\mbox{if $x\in D_N$}	\nonumber
\end{eqnarray}
with $v_0(x)\le u_0(x)$.  For $x\in\mathbf{Z}^d\setminus D_N$, let $v(t,x)=0$. 
Then, with probability 1, for all $(t,x)\in [0,\infty)\times\mathbf{Z}^d$ we
have $v(t,x)\le u(t,x)$.  
\end{lemma}

\noindent\textbf{Proof.  }
The lemma follows from standard comparison arguments.  See, for example,
\cite{kot92}.  \qed

As in Lemma \ref{l.int-eq}, the solution $v$ to (\ref{cutoff}) will satisfy
\begin{eqnarray}
\label{vint-eq}
v(t,x)&=&\sum_{y\in\mathbf{Z}^d}G_N(t,x-y)v_0(y) \\
	&&+\int_0^t\sum_{y\in\mathbf{Z}^d}G_N(t-s,x-y)
      \cdot v^\gamma(s,y)dB_y(s),\ x\in D_N, \nonumber 
\end{eqnarray}
where $G_N(t,x)$ is the fundamental solution of the discrete heat equation 
on  $\mathbf{Z}^d$ with $0$ boundary conditions on $\partial D_N$.

We will at times use the following consequence of Jensen's inequality. 
Let $M>0$, suppose that $p>1$, and that $a_1,\dots,a_M$ are non-negative
real numbers.  Then,
\begin{eqnarray}
\label{jensen}
\sum_{k=1}^M a_k^p&=&M\left(\sum_{k=1}^M a_k^p\frac{1}{M}\right)	\\
&\ge&M\left(\sum_{k=1}^M a_k\frac{1}{M}\right)^p	\nonumber	\\
&=&M^{1-p}\left(\sum_{k=1}^M a_k\right)^p.		\nonumber
\end{eqnarray}

For the following lemma, let $\# S$ denote the cardinality of the set $S$.
\begin{lemma}
\label{vdecay}
Suppose that $v(t,x)$ satisfies (\ref{cutoff}), and let 
\[
V(t)=\sum_{x\in D_N}v(t,x).
\]
Let $\partial^- D_N$ be the points in $D_N$ which are nearest neighbors 
to $\partial D_N$.
Then there exists a Brownian motion $B(t)$ and a predictable functional
$H(t)\ge (3N)^{-(2\gamma-1)d/2}V^\gamma(t)$ 
such that $V(t)$ satisfies the following stochastic differential
equation:
\begin{eqnarray}
\label{mass}
dV(t)&=&-\sum_{x\in\partial^- D_N}v(t,x)[2d-\#(\mathcal{N}(x)\cap D_N)]
dt\\
& &\qquad\qquad+H(t)dB(t)	\nonumber\\
V(0)&=&\sum_{x\in D_N}v(0,x)	\nonumber
\end{eqnarray}
\end{lemma}

\noindent\textbf{Proof.  }
Summing over $x\in D_N$, combining the Brownian motions and using 
(\ref{jensen}) with $p=2\gamma$, we get the
lemma.  \qed

The following lemma is a special case of Lemma 3.4 of \cite{mp92}.
\begin{lemma}
\label{die-out}
Let $A>0$, $\gamma\in (1/2,1)$, and let $Z(t)$ satisfy
\begin{eqnarray*}
dZ(t)&=&AZ^\gamma(t)dB(t)	\\
Z(0)&=&z_0>0.
\end{eqnarray*}
Let $\tau$ be the first time $t$ that $Z(t)=0$, and let $\tau=\infty$ if
$Z(t)$ never reaches 0.  There is a constant $C(\gamma)$ depending only
on $\gamma$ such that  
\[ P(Z(t)>0)=
P(\tau>t)\le C(\gamma)A^{-2}z_0^{2-2\gamma}t^{-1}.
\]
\end{lemma}
Lemma \ref{die-out} has the following simple corollary.

\begin{lemma}
\label{a3} 
Let $X_t$, $t\le T$,  be a 
continuous non-negative supermartingale with martingale part 
\[
\int_{0}^{t}H_sdB_s
\]
for some Brownian motion $B_t$ and for some predictable process $H_t$. 
If $L\ge 1$ and $A>0$, then 
\[
P\left( X_{T}>0,\int_0^T1(H_t\ge AX_t^\gamma)dt\ge T/L
\right)\le C(\gamma,L)A^{-2}X_0^{2-2\gamma}T^{-1}.
\]
\end{lemma}

\noindent\textbf{Proof.  }
Let $\sigma (t)=\int_0^{t\wedge T_X} H_u^2A^{-2}X_u^{-2\gamma}du$, where
$T_X=\inf\{t:X_t=0\}$, and define $\tau(s)=\inf\{t:\sigma(t)>s\}$, if 
$s<\sigma(T_X)$;  and $\tau(s)=T_X$, if $s\ge\sigma(T_X)$.  Then 
$Z_t=X_{\tau(t)}$ is a right continuous non-negative supermartingale with 
continuous martingale part given by $\int_0^{t\wedge T_Z}
AZ_s^\gamma\tilde dB_s$ for some Brownian motion $\tilde B$.  $Z$ may only 
have negative jumps.  Let $Y_t$ be the unique solution of 
\[ 
Y_t=X_0+\int_0^tAY_s^\gamma d\tilde B_s.
\]
A well-known comparison theorem (Theorem V.43.1 of \cite{rw87} is easily 
modified to cover this case) shows that $Z_{t-}\le Y_t$ for all $t\ge 0$ a.s.
and so it follows easily that $X_t\le Y_{\sigma(t)}$ for all $t\ge 0$ a.s.  
Clearly 
\[
\int_0^T1(H_t\ge AX_t^\gamma)dt\ge T/L
\]
and $T<T_X$ imply
$\sigma(T)\ge T/L$.  Recalling that both $X$ and $Y$ will stick at 0 after 
they first hit $0$, we see that
\begin{eqnarray*}
P\left(X_T>0, \int_0^T1(H_t\ge AX_t^\gamma)dt\ge T/L\right)
&\le & P\left(\sigma (T)\ge T/L, Y_{\sigma(T)}>0\right)\\
&\le & P\left(Y_{T/L}>0\right), 
\end{eqnarray*}
and Lemma \ref{die-out} completes the proof. \qed

\begin{lemma}
\label{even}
Suppose that $f$ and $g$ are non-negative functions 
on $\mathbf{Z}^d$ such that 
\begin{eqnarray*}
\sum_x f(x)&\le& M,   \\
\sum_x g(x)f(x)&\le& K.
\end{eqnarray*}  
Then
\[
\sum_x g(x)^{2-2\gamma}f(x)\le K^{2-2\gamma}M^{2\gamma-1}.
\]
\end{lemma}

\noindent\textbf{Proof.  }
By H\"older's inequality,
\begin{eqnarray*}
\sum_x g(x)^{2-2\gamma}f(x)&=& \sum_x \left[\left(g(x)f(x)\right)^{2-2\gamma}
		f(x)^{1-(2-2\gamma)}\right]  \\
   &\le& \left[\sum_x g(x)f(x)\right]^{2-2\gamma}
	\left[\sum_x f(x)\right]^{2\gamma-1}  \\
   &\le& K^{2-2\gamma}M^{(2\gamma-1)}.\qed
\end{eqnarray*}

Our final result shows that we can split up solutions at $t=0$ in an
appropriate manner.  First we observe that since $2\gamma>1$, if $a,b>0$
then 
\[
(a+b)^{2\gamma}\ge a^{2\gamma}+b^{2\gamma}.
\]

\begin{lemma}
\label{decomp2}
Fix $n\ge 0$, and let $u_0:\mathbf{Z}^d\to[0,\infty)$. 
For each 
$1\le i\le n$, let $S_i$ be a finite subset of $\mathbf{Z}^d$, and let 
$w_{i,0}(\cdot)$ be a non-negative function supported on $S_i$.  Assume that
\[
\sum_{i=1}^nw_{i,0}(x)\le u_0(x).
\]
Then on some filtered probability space, $(\Omega,\mathcal{F},\mathcal{F}_t,P)$, 
we may define independent $\mathcal{F}_t$-Brownian motions 
$\{B_{i,x}:1\le i\le n, x\in \mathbf{Z}^d\}$, independent 
$\mathcal{F}_t$-Brownian motions 
$\{B_{x}:x\in \mathbf{Z}^d\}$ (the 2 collections will not be mutually 
independent), and non-negative $(\mathcal{F}_t)$-predictable continuous 
functions $u(\cdot,x)$ and $w_i(\cdot,x)$, $H_i(\cdot,x)$ 
$(i\le n, x\in\mathbf{Z}^d)$ such that the following holds.  First, the 
$w_i(t,x)$ satisfy
\begin{eqnarray}
\label{h.equ}
dw_i(t,x)&=&\Delta w_i(t,x)dt+H_i(t,x)dB_{i,x}(t)\quad\mbox{if $x\in S_i$}\\
w_i(t,x)&=&0 \quad\mbox{ if $x\not\in S_i$}	\nonumber\\
w_i(0,x)&=&w_{i,0}(x)		\nonumber
\end{eqnarray}
and $u(t,x)$ satisfies (\ref{system}).  Secondly, 
\[
H_i(t,x)\ge w_i(t,x)^\gamma.
\]
Finally, with probability 1, for all $t\ge 0$, $x\in\mathbf{Z}^d$, we have
\[
\sum_{i=1}^nw_i(t,x)\le u(t,x).
\]
\end{lemma}

\noindent\textbf{Proof.  }
We give a brief outline of the proof, leaving the details to the reader.
Observe that if the continuous function 
$h_i:\mathbf {R}^n\to \mathbf{R}_+$ is defined by 
\begin{equation}
\label{Hdef}
h_i(w_1,{\ldots}, w_n)
     =\frac{w_i^\gamma\left(\sum_{k=1}^{n}w_k\right)^\gamma}
	{\left(\sum_{k=1}^{n}w_k^{2\gamma}\right)^{1/2}}
	\mathbf{1}\left(\sum_{k=1}^{n}w_k> 0\right),
\end{equation}
then
\[
h_i(w_1,{\ldots}, w_n)\ge w_i^\gamma
\]
and
\[
\sum_{i=1}^{n}h_i^2(w_1,{\ldots} ,w_n)
      =\left(\sum_{i=1}^{n}w_i\right)^{2\gamma}.
\]
Let $h_{i}^{(m)}$ be a sequence of Lipschitz functions on $\mathbf {R}^n$, 
which converge to $h_{i}$ 
uniformly.  More specifically, we define a function 
$\phi_m:  \mathbf{R}_+\to\mathbf{R}$ as follows.  Let $\phi_m(w)$ agree with 
$w^\gamma$ on $[1/m,\infty)$ and at $0$, and define $\phi_m(w)$ by linear
interpolation on $(0,1/m)$.  Define $h_i^{(m)}$ as $h_i$ but with 
$\phi_m(w_i)\phi_m (\sum_{k=1}^n w_k)$ in the numerator.  It follows 
that $\sum_{i=1}^{n}{h_i^{(m)}}^2(w_1,{\ldots}, w_n)$ converges to 
$(\sum_{i=1}^{n}w_i)^{2\gamma}$ uniformly on compacts.  Let $w_i^{(m)}(t,x)$ 
be the unique solution of
\begin{eqnarray*}
dw_i^{(m)}(t,x)&=&\Delta w_i^{(m)}(t,x)dt  
        +h_i^{(m)}(w_1^{(m)}(t,x),{\ldots} ,w_n^{(m)}(t,x))dB_{i,x}(t) \\
&&      \qquad\qquad\qquad\mbox{if $x\in S_i$}\\
w_i^{(m)}(t,x)&=&0 \quad\mbox{ if $x\not\in S_i$}	\\
w_i^{(m)}(0,x)&=&w_{i,0}(x).		
\end{eqnarray*}
Let $z_i^{(m)}(t,x)$ be the unique solution of
\begin{eqnarray*}
dz_i^{(m)}(t,x)&=&\Delta z_i^{(m)}(t,x)dt
     +h_i^{(m)}(z_1^{(m)}(t,x),{\ldots} ,z_n^{(m)}(t,x))dB_{i,x}(t) \\
z_i^{(m)}(0,x)&=&\bar w_{i,0}(x),		
\end{eqnarray*}
where $\bar w_{i,0}(x)\ge w_{i,0}(x)$ are chosen such that 
\[
\sum_{i=1}^{m}\bar w_{i,0}(x)=u_0(x).
\]
Then, standard comparison theorems show that with probability $1$,
\[
0\le w_i^{(m)}(t,x)\le z_i^{(m)}(t,x)\hbox{ for all }t,x.  
\]
Furthermore, we can
take a subsequence $m_k$ such that $(w_i^{(m_k)}(t,x),z_i^{(m_k)}(t,x))$
converges in distribution to $(w_i(t,x),z_i(t,x))$, where $(w_i)$ satisfy
(\ref{h.equ}) with $H_i(t,x)=h_i(w_1(t,x),{\ldots}, w_n(t,x))$ and 
$\sum_{i=1}^{n}z_i(t,x)$ satisfies (\ref{system}).  Clearly with
probability 1, 
\begin{equation}
0\le w_i(t,x)\le z_i(t,x) \hbox{ for all }t,x,i.\nonumber
\end{equation}
This implies
Lemma \ref{decomp2}.
\qed

\section{Proof of Theorem \ref{t1}}
  \setcounter{equation}{0}

Our proof depends heavily on the decomposition in Lemma \ref{decomp2}.
We emphasize that this decomposition is a weak existence theorem, so we
are always dealing with a changing set of Brownian motions and solutions.
  
First, we set up some notation which we will use in the proof.  Fix
$\varepsilon>0$.
Let $\mathcal{K}$ be the event that for some $T>0$, $u(T,x)=0$ for all 
$x\in\mathbf{Z}^d$.  For $T\ge 0$, let $\mathcal{K}(T)$ be the event that 
$u(T,x)=0$ for all $x\in \mathbf{Z}^d$. Note that as $T\uparrow \infty$, the
events $\mathcal{K}(T)$ increase to $\mathcal{K}$. In fact, we show that there
exists a time $t_\infty$ such that 
\begin{equation}
\label{main-est}
P\left( \mathcal{K}(t_\infty)^c\right)<\varepsilon.  
\end{equation}
Then (\ref{main-est}) implies Theorem 1.  From now on, we fix
$\varepsilon>0$ and concentrate on proving (\ref{main-est}).  

Let $\ell\in \mathbf {N}$ 
be so large
that 
\begin{equation}
\label{ell2}
d\cdot\gamma^{\ell-1}<\frac{1}{8},
\end{equation}
and then choose $0<\delta$ small enough so that
\begin{equation}
\label{delta}
\delta<{1\over 2}\gamma^{\ell-1}(1-\gamma).
\end{equation}
We will specify an integer $n_0>0$ later.  Let 
\[
N_n=2^n.
\]
We define $0=t_{n_0}<t_{n_0+1}<\ldots$ inductively as follows.  
Let $\bar m$ be the smallest integer $m$ such that 
\[
m\ell\ge 2^{n_0/2}
\]
and let 
\[
t_{n_0+1}={\bar m}\ell.
\]
If $n>n_0$ and if we are given $t_n$, we let
\[
t_{n+1}=t_n+2^{-n/2}.
\]
Clearly, there exists a finite accumulation point
\[
t_\infty=\lim_{n\to\infty}t_n<\infty.
\]

Fix $K>0$ and let $\mathcal{A}_0=\mathcal{A}_0(K)$ be the event that 
\[
\sup_{t\ge 0}\sum_{x\in\mathbf{Z}^d}u(t,x)\le K.
\]
Let $\tilde M_t=\sum_{x\in\mathbf{Z}^d}u(t,x)$.  Note that the integrability
in (\ref{rw.0}) easily gives that $\tilde M_t$ is a continuous non-negative
martingale.  Let ${\tilde\tau}={\tilde\tau}(T,K)$ be the first time 
$t\le T$ that $\tilde M_t\ge K$.  If there is no such time, let 
${\tilde\tau}=T$.  Using the optional sampling theorem and Markov's
inequality, we get
\begin{eqnarray}
\label{Kbound}
P\left(\mathcal{A}_0(K)^c\right)
	&=&\lim_{T\to\infty}P\left(\tilde M_{{\tilde\tau}(T,K)}\ge K\right)	\\
	&\le& \frac{\tilde M_0}{K}.	\nonumber
\end{eqnarray}

For $n\ge n_0$, $t_{n}\le t\le t_{n+1}$, 
$x\in\mathbf{Z}^d$, we will inductively
define a sequence of random functions $v_n(t,x)\le u(t,x)$ as follows. 
Since the definition of $v_n(t,x)$ for $n>n_0$ is the simplest, we start
with that case.  Suppose that we have defined $v_{n-1}(t,x)\le u(t,x)$ for 
$t_{n-1}\le t\le t_{n}$.  Let 
\begin{equation}
\label{defvn}
v_{n}(t_n,x)=u(t_n,x)-v_{n-1}(t_n,x). 
\end{equation}
For $t_n<t\le t_{n+1}$, let $v_n(t,x)$ satisfy equation
(\ref{cutoff}) (of Lemma \ref{decomp1}), with $N=N_n$, so that 
$v_n(t,x)\le u(t,x)$ by Lemma \ref{decomp1}.

Now we give the more complicated definition of $v_{n_0}(t,x)$.  For 
$1\le m\le\bar m$, we call the time intervals 
$[(m-1)(\ell),m(\ell))$ stages, and we call the subintervals 
$[k,k+1)\subset [(m-1)(\ell),m(\ell))$ substages.
In order to define $v_{n_0}(t,x)$, we first define a collection of
functions $w_{k,z}(t,x)$ for $k\ell\le t\le (k+1)\ell$, $0\le k<{\bar{m}}$, 
$z\in D_{N_{n_0}}$, satisfying
\begin{equation}
\label{wsum} \sum_{z\in D_{N_{n_0}}} w_{k,z}(t,x)\le u(t,x), \quad k\ell \le t
\le (k+1)\ell.
\end{equation}

We let $x+D_k$ denote the set $\{x+y: y\in D_k\}$.
Let 
\[
\bar{N}_n=2^{\delta n}.
\]
Assume that either $k=0$ or that $w_{k-1,z}$ has already been defined
for $z\in D_{N_{n_0}}$ and satisfies (\ref{wsum}).  Let
\[
w_{0,z}(0,x)=u(0,z)\mathbf{1}(x=z),
\]
and let 
\[
w_{k,z}(k\ell,x)
	=\sum_{y\in D_{N_{n_0}}}w_{k-1,y}(k\ell,z)\mathbf{1}(x=z).
\]
Therefore (\ref{wsum}) holds for $t=k\ell$.
In either case, using Lemma \ref{decomp2} with $S_z=z+D_{\bar N_{n_0}}$, 
$z\in D_{N_{n_0}}$,  
we let $w_{k,z}(t,x)$ satisfy
the following equation for 
$k\ell\le t\le (k+1)\ell$.
\begin{eqnarray*}
dw_{k,z}(t,x)&=&\Delta w_{k,z}(t,x)dt+H_z(t,x)dB_{z,x}(t)  \quad
x\in z+D_{\bar{N}_{n_0}}\\
w_{k,z}(t,x)&=&0,\qquad x\notin z+ D_{\bar{N}_{n_0}}, 
\end{eqnarray*}
where $H_z(t,x)$ is as in Lemma \ref{decomp2}, and 
$\{B_{z,x}(t)\}_{z,x}$ are independent Brownian motions.
Lemma \ref{decomp2} now implies (\ref{wsum}) and our inductive construction
is complete.

Finally, for $0\le t \le t_{n_0+1}$, we define $v_{n_0}(t,x)$ as follows. 
Let $D^0_{N_{n_0}}$ be the set of those points $z\in D_{N_{n_0}}$ 
such that 
\[
z+D_{\bar{N}_{n_0}}\subset D_{N_{n_0}}.
\]
The reason for defining $D^0_{N_{n_0}}$ is that we do not want to include 
any points $z$ for which $w_{k,z}(t,x)$ has any support outside of 
$D_{N_{n_0}}$.  If $x\not\in D^0_{N_{n_0}}$, let $v_{n_0}(t,x)=0$.  
If $x\in D^0_{N_{n_0}}$ and $t\le {\bar m}\ell=t_{n_0+1}$, 
choose $0\le k<\bar {m}$ such that 
$k\ell\le t\le (k+1)\ell$ and let 
\[
v_{n_0}(t,x)=\sum_{z\in D^0_{N_{n_0}}}w_{k,z}(t,x)\le u(t,x).
\]
Extend $v_n$, $w_{k,z}$ to be identically $0$ outside their initial 
domains of definition and let $\mathcal{F}_t$ 
be the right-continuous filtration
generated by the processes $v_n$, $w_{k,z}$, $u$, and $B_{z,x}$ up to time
$t$ as $n, k, z, x$ range through their respective domains.  
Now we label the mass that has ``leaked out".  For $n\ge n_0$, let 
\[
M_n=\sum_{x\in\mathbf{Z}^d}[u(t_{n+1},x)-v_n(t_{n+1},x)].
\]

We will often work with the following sets.
\begin{definition}
For $n\ge n_0$ let $\mathcal{A}_{1,n}$ be the event that 
\[
v_n(t_{n+1},x)=0 \hbox{ for all }x\in\mathbf{Z}^d,
\]
and let $\mathcal{A}_{2,n}$ be the event that 
\[
M_{n}\le 2^{-2^{\delta n}}.
\] 
Define $\mathcal{A}_{i,n_0-1}$ to be 
the entire space for $i=1$ or $2$. 
\end{definition}
Our definitions imply that 
\[
\mathcal{K}(t_\infty)\supset\mathcal{A}_0(K)
	\cap\left[\bigcap_{n=n_0}^\infty \mathcal{A}_{1,n}\right]
	\cap\left[\bigcap_{n=n_0}^\infty \mathcal{A}_{2,n}\right].
\]

The following lemma plays an essential role in the proof of Theorem 1.

\begin{lemma}
\label{l.est-a}
If $n_0$ is large enough and $n\ge n_0$, then
\[ 
P\left(\mathcal{A}_{2,n}^c\cap\mathcal{A}_0(K)\cap\mathcal{A}_{1,n-1}\right)
	\le 2^{-2^{\delta n}},
\] 
and
in particular, for $n_0$ large enough,
\begin{eqnarray}
\label{est-a}
P\left(\bigcup_{n=n_0}^\infty\left[\mathcal{A}_{2,n}^c
	\cap\mathcal{A}_0(K)\cap\mathcal{A}_{1,n-1}\right]\right)
	 &\le& \sum_{n=n_0}^\infty 2^{-2^{\delta n}} \\
 	&\le& \frac{\varepsilon}{4}.  \nonumber
\end{eqnarray}
\end{lemma}

\noindent
\textbf{Proof.  }
To begin the proof of Lemma \ref{l.est-a}, we consider the case 
$n>n_0$.  The key to the proof is the observation that on the set 
$\mathcal{A}_{1,n-1}$ we have $v_n(t_n,\cdot)=u(t_n,\cdot)$.  This follows
from the definitions of $v_n$ and the event $\mathcal{A}_{1,n-1}$. 
Let $\xi_t^x$ be our original continuous time random walk, started from $x$,
and let $\tau_n^x$ be the first time $t$ that $\xi_t^x\not\in D_{N_n}$.
The integral equations (\ref{int-eq}) and (\ref{vint-eq})
imply that
on $\mathcal{A}_{1,n-1}\in\mathcal{F}_{t_n}$,
\begin{equation}
\label{Mn}
E\left( \sum_xu(t_{n+1},x)-v_n(t_{n+1},x)\bigg|\mathcal{F}_{t_n}\right)=\sum_x
u(t_n,x)P(\tau_n^x<2^{-n/2}).
\end{equation}
Now use (\ref{int-eq}) again to see that
\begin{eqnarray*}
\lefteqn{E\left(M_n\mathbf{1}(\mathcal{A}_0(K)\cap\mathcal{A}_{1,n-1})\right)}
\\
&\le &  \mathbf{1}(\sum_x u_0(x)\le K)E\left(\mathbf{1}({\mathcal{A}_{1,n-1}})
E\left(\sum_xu(t_{n+1},x)-v_n(t_{n+1},x)\bigg|\mathcal{F}_{t_n}\right)\right)\\
&\le &
\mathbf{1}(\sum_xu_0(x)\le K)E\left(\sum_x u(t_n,x)P(\tau_n^x<2^{-n/2})
\right)\\
&\le & 
\mathbf{1}(\sum_xu_0(x)\le K)\sum_xu_0(x)P(\tau_n^x<t_n+2^{-n/2}=t_{n+1})\\
&\le & K\sup_{x\in\mathbf{F}}P(\tau_n^x<t_{n+1}). 
\end{eqnarray*}
Denote
\[
\|\mathbf{F}\|=\max_{x\in\mathbf{F}}|x|.    
\]
Let $S_t$ be the number of steps that $\xi_s$ has taken for $s\le t$.  Of
course, since the steps are of size 1, we have that 
\[
P\left(\tau_n^x\le t_{n+1}\right)\le P\left(S_{t_{n+1}}>2^n-\|\mathbf{F}\|\right).
\]
Recall the definition of $M_n$.  Also, from the definition of $t_{n_0}$
and  $t_n$, we see that for $n_0\ge n(\ell)$,
\[
t_{n+1}\le t_{n_0}+\sum_{m=n_0}^{\infty}2^{-m/2}\le 2^{1+n_0/2}.
\]
Using (\ref{poisson}), we conclude that
\begin{eqnarray}
\label{m-est}
E\left[M_{n}\mathbf{1}(\mathcal{A}_0(K)\cap\mathcal{A}_{1,n-1}\right)]
	&\le& KP\left(S_{t_{n+1}}>2^n-\|\mathbf{F}\|\right)	\\
	&\le& \frac{CK(2dt_{n+1})^{2^n-\|\mathbf{F}\|}}
		{\sqrt{2^n-\|\mathbf{F}\|}\left(2^n-\|\mathbf{F}\|\right)
		  ^{2^n-\|\mathbf{F}\|}e^{-2^n+\|\mathbf{F}\|}} 
		\nonumber\\
	&\le& 8^{-2^n}	\nonumber
\end{eqnarray}
if $n$ is large enough, and thus if $n_0$ is large enough.
Using Markov's inequality, we have
\[
P\left(\{M_{n}>2^{-2^n}\}\cap\mathcal{A}_0(K)\cap\mathcal{A}_{1,n-1}\right)\le 
4^{-2^n}.
\]
This proves Lemma \ref{l.est-a} for $n>n_0$ because $\delta<1$.

Next, we turn to the case $n=n_0$. Note that $M_{n_0}$ consists of 2
kinds of mass.  Let $M^{'}_{n_0}$ refer to the first kind of mass,
which escapes from each of the small cubes $x+D_{{\bar N}_{n_0}}$.
Let $M^{''}_{n_0}$ refer to the second kind of mass, which escapes 
from the large cube $D_{N_{n_0}}$ or becomes part of the functions
$w_{k,z}(t,x)$, for those $z\in D^0_{N_{n_0}}$.  To be precise, let 
\[
M^{'}_{n_0}=\sum_{x\in\mathbf{Z}^d}\left[
u(t_{n_0+1},x)-\sum_{z\in D_{N_{n_0}}}w_{{\bar m}-1,z}(t_{n_0+1},x)\right]
\]
and let 
\[
M^{''}_{n_0}
=\sum_{x\in\mathbf{Z}^d}\left[\sum_{z\in D_{N_{n_0}}-D^0_{N_{n_0}}}
	w_{{\bar m}-1,z}(t_{n_0+1},x)\right].
\]
Clearly $M_{n_0}=M^{'}_{n_0}+M^{''}_{n_0}$.

First we deal with $M^{'}_{n_0}$.  Let $\bar{\tau}$ be the first time
$t<\ell$ that $\xi_t\notin D_{\bar{N}_{n_0}}$ where $\xi_0=0$.
If there is no such
time, let $t=\ell$.  
Using the analogue of (\ref{vint-eq}) for the $w_{k,z}$, as in (\ref{Mn}),
we get
\begin{eqnarray*}
\lefteqn{E\left(\sum_x\left[u(t_{n_0+1},x)-\sum_zw_{{\bar m}-1,z}
      (t_{n_0+1},x)\right]\Bigg|
\mathcal {F}_{({\bar m}-1)\ell}\right)}\\
&= &\sum_x\left[u(({\bar m}-1)\ell,x)-\sum_zw_{{\bar m}-1,z}(({\bar m}-1)\ell
,x)P({\bar \tau}\ge \ell)\right]\\
&= &\sum_x\left[u(({\bar m}-1)\ell,x)-\sum_zw_{{\bar m}-1,z}(({\bar m}-1)\ell,
x)\right]\\
&\ & \qquad+\sum_{x,z}w_{{\bar m}-1,z}(({\bar m}-1)\ell,x)P(\bar \tau<\ell)
 \\
&\le & \sum_x\left[u(({\bar m}-1)\ell,x)-\sum_zw_{{\bar m}-2,z}(({\bar m}-1)
\ell,x)\right]\\
&\ &\qquad +\sum_xu(({\bar m}-1)\ell,x)P(\bar\tau<\ell). 
\end{eqnarray*}
In the last line we have used (\ref{wsum}) and the fact that at times $k\ell$,
the redistribution of the mass among the $w_{k,z}$'s preserves the total
mass of $\sum_zw_{k,z}$.  Therefore
\begin{eqnarray*}
\lefteqn{E\left(\sum_x\left[u(t_{n_0+1},x)-\sum_zw_{{\bar m}-1,z}(t_{n_0+1},x)\right]
\mathbf{1}(\mathcal{A}_0(K))\right)} \\
&\le& E\left( \sum_x\left[u(({\bar m}-1)\ell,x)-\sum_zw_{{\bar m}-2,z}(({\bar m
}-1)\ell,x)\right]\mathbf{1}(\mathcal{A}_0(K))\right)\\
&\ &\qquad +E\left(\sum_xu(({\bar m}-1)\ell,x)\mathbf{1}
(\sum_x \bar u((\bar m-1)\ell,x)\le K)\right)P(\bar\tau<\ell).
\end{eqnarray*}
The last term is at most $KP(\bar \tau<\ell)$.  Now iterate the above 
$\bar m$ times, noting that $\sum_xu(0,x)=\sum_x\sum_zw_{0,z}(0,x)$,
and argue as in (\ref{m-est}) using (\ref{poisson}) 
to get (for $n_0$ large again)
\begin{eqnarray*}
E\left[M^{'}_{n_0}\mathbf{1}(\mathcal{A}_0(K))\right]
  &\le& K\bar m P(\bar \tau <\ell)\\
&\le& K2^{n_0/2}P(S_\ell\ge 2^{\delta n_0})\\
	&\le& CK2^{n_0/2}
	    \frac{(2d\ell)^{2^{\delta n_0}}}
		{\sqrt{2^{\delta n_0}}\left(2^{\delta n_0}\right)
		^{2^{\delta n_0}}e^{-2^{\delta n_0}}}	\\
	&\le& 8^{-2^{\delta n_0}},
\end{eqnarray*}
if $n_0$ is large enough.  
Again using Markov's inequality, we have 
\begin{equation}
\label{mm1}
P\left(M^{'}_{n_0}>\frac{1}{2}2^{-2^{\delta n_0}},\mathcal{A}_0(K)\right)
	\le \frac{1}{2}2^{-2^{\delta n_0}}.
\end{equation}

Now we consider $M^{''}_{n_0}$.  
On the last interval $[(\bar m-1)\ell,\bar m\ell]$, 
\[
\sum_{z\in D_{N_{n_0}}\setminus D^0_{N_{n_0}}}
\sum_x w_{\bar m-1,z}(t,x)
\]
is a supermartingale and so
\begin{eqnarray*}
\lefteqn{E\left(M^{''}_{n_0}\mathbf{1}(\mathcal{A}_0(K))\right)}\\
&\le&E\left(\sum_{z\in D_{N_{n_0}}-D_{N_{n_0}}^0}\sum_xw_{\bar m-1,z}
((\bar m-1)\ell,x)\mathbf{1}(\langle u_0,\mathbf{1}\rangle \le K)\right)\\
&=&E\left(\sum_{z\in D_{N_{n_0}}-D_{N_{n_0}}^0}\sum_x\sum_{y\in D_{N_{n_0}}}
w_{\bar m-2,y}((\bar m-1)\ell,x)\mathbf{1}(x=z)\right)  \\
  && \qquad\cdot\mathbf{1}(\langle u_0,\mathbf{1}\rangle \le K)\\
&=& E\left(\sum_{x\in D_{N_{n_0}}-D_{N_{n_0}}^0}\sum_{y\in D_{N_{n_0}}}
w_{\bar m-2,y}((\bar m-1)\ell,x)\right)\mathbf{1}(\langle u_0,\mathbf{1}\rangle
\le K).
\end{eqnarray*}
Use (\ref{wsum}) to bound the above by ($n_0$ large enough)
\begin{eqnarray*}
\lefteqn{
E\left(\sum_{x\in D_{N_{n_0}}-D_{N_{n_0}}^0}u((\bar m-1)\ell,x)\right)
\mathbf{1}(\langle u_0,\mathbf{1}\rangle \le K)}\\
&\le & K\sup_{x\in \mathbf{F}}
P(\xi_x \hbox{ exits }D^0_{N_{n_0}} \hbox{ before time }(\bar m -1)\ell)\\
&\le& KP(S_{2^{n_0/2}}>2^{n_0-1}).
\end{eqnarray*}
Another application of (\ref{poisson}) (as in (\ref{m-est}))
shows that for $n_0$ large enough 
the above is at most 
$8^{-2^{n_0}}$, and therefore, by Markov's inequality,
\begin{equation}
\label{mm2}
P\left(M^{''}_{n_0}>\frac{1}{2}2^{-2^{\delta n_0}},\mathcal{A}_0(K)\right)
	\le \frac{1}{2}2^{-2^{\delta n_0}}.
\end{equation}
Putting together (\ref{mm1}) and (\ref{mm2}), we obtain
\begin{eqnarray*}
P\left(M_{n_0}>2^{-2^{\delta n_0}},\mathcal{A}_0(K)\right)
 &\le&P\left(M^{'}_{n_0}>\frac{1}{2}2^{-2^{\delta n_0}},\mathcal{A}_0(K)\right)\\
 &&+P\left(M^{''}_{n_0}>\frac{1}{2}2^{-2^{\delta n_0}},\mathcal{A}_0(K)\right) \\
 &\le& 2^{-2^{\delta n_0}}.
\end{eqnarray*}
This proves Lemma \ref{l.est-a}.  \qed

Our next goal is to estimate the probability of $\mathcal{A}_{1,n}^c$.  
Let 
\[
V_n(t)=\sum_{x\in D_{N_n}}v_n(t_n+t,x).
\]

\begin{lemma}
\label{a2}
If $n_0>n(\delta,\varepsilon)$, then for $n>n_0$,
\[
P\left( \mathcal{A}_{1,n}^c\cap\mathcal{A}_{2,n-1}\right) 
	\le 2^{n_0-n}\frac{\varepsilon}{4}.
\]
\end{lemma}

\noindent\textbf{Proof.  }
Our argument uses Lemma \ref{a3}.  We need to bound
\[
P\left(\mathcal{A}_{1,n}^c\cap\mathcal{A}_{2,n-1}\right)
	=P\left(V_n(2^{-n/2})>0,\mathcal{A}_{2,n-1}\right).
\]
By Lemma \ref{vdecay}, on the event $\mathcal{A}_{2,n-1}$, we can write
\begin{eqnarray*}
dV_n(t)&\le&H(t)dB(t),\qquad 0<t<t_{n+1}-t_n=2^{-n/2}	\\
V_n(0)&\le&2^{-2^{\delta (n-1)}},
\end{eqnarray*}
where
\begin{eqnarray*}
H(t)&\ge&(3N_n)^{-(2\gamma-1)d/2}V_n(t)^\gamma\\
   &=& C2^{-nd(\gamma-1/2)}V_n(t)^\gamma.
\end{eqnarray*}
Then, by Lemma \ref{a3}, we have
\[
P\left( V_n(2^{-n/2})>0,\mathcal{A}_{2,n-1}\right)
  \le C2^{nd(2\gamma-1)}2^{-2^{\delta (n-1)}(2-2\gamma)}2^{n/2}.
\]
This proves Lemma \ref{a2}, if $n_0$ is large enough.  \qed

Next, we treat the more complicated case of $n=n_0$.
\begin{lemma}
\label{n0}
If $n_0\ge n(\varepsilon,K)$, then
\[
P\left( \mathcal{A}_{1,n_0}^c\cap\mathcal{A}_0(K)\right) 
    \le \frac{\varepsilon}{4}.
\]
\end{lemma}

To prove this we will deal with the stages 
(of length $\ell$) and the substages 
(of length 1) which we defined earlier.  We first show that 
for at least half of the stages, 
at the
end of the $\ell-1$ substages, there are only a small number of sites 
$z\in D_{N_{n_0}}$ such that $w_{k,z}$ is still alive.  (Recall that we say
a function is alive if it is not identically 0).  
To state this key lemma precisely, for $0\le k \le \bar m$, we let
\[
W_{k,z}(t)=\sum_{x\in z+D_{{\bar N}_{n_0}}}w_{k,z}(t,x),
\quad k\ell\le t\le (k+1)\ell,
\]  
and for $0\le j<\ell$ and $k$ as above, set 
\[
\mathcal{A}_{3,k}(j)=\left\{\sum_{z\in D_{N_{n_0}}}\mathbf{1}(W_{k,z}(k\ell+j)>0)
\le 3^{d\gamma^j}2^{n_0d\gamma^j}\right\}.\nonumber
\]
Then we have the following.

\begin{lemma}
\label{decrease}
If $n_0\ge n(\ell,K,\gamma,\varepsilon)$ then
\[
P\left(\left\{\sum_{k=0}^{\bar m -1}\mathbf{1}(\mathcal{A}_{3,k}(\ell-1)^c)\ge \bar m/2\right\}
\cap\mathcal{A}_0(K)\right)<\varepsilon/8 \nonumber.
\]
\end{lemma}

Assume for the moment that Lemma \ref{decrease} holds and let us give 
the 
\newline
\noindent\textbf{Proof of Lemma \ref{n0}.  }
Note that, as in Lemma \ref{vdecay}, $W_{k,z}$ is a non-negative
supermartingale with martingale part $H_{k,z}dB$, where 
\begin{eqnarray}
\label{lbH}
H_{k,z}(t)&\ge&\sum_{x\in z+D_{\bar N_{n_0}}}w_{k,z}(t,x)^\gamma \\
\nonumber 
&\ge&\left(\sum_{x\in z+D_{\bar N_{n_0}}}
w_{k,z}(t,x)^{2\gamma}\right)^{1/2} \\ \nonumber
 &\ge& \left(3^d2^{\delta n_0d}\left(\frac{W_{k,z}(t)}
{3^d2^{\delta n_0d}}\right)
	^{2\gamma}\right)^{1/2} \qquad \hbox{(by (\ref{jensen}))}\\
 &=& C2^{-\delta n_0d(\gamma-1/2)}W_{k,z}(t)^\gamma. \nonumber
\end{eqnarray}
Lemma \ref{vdecay} also shows that $V_{n_0}(t)=\sum_{z\in D_{N_{n_0}}}
W_{k,z}(t)$ (for $k\ell\le t <(k+1)\ell$) is a continuous supermartingale
with martingale part $H_tdB_t$, where by (\ref{lbH}) and Jensen's inequality, 
for $k\ell\le t <(k+1)\ell$
\begin{eqnarray*}
H_t^2=\sum_zH_{k,z}(t)^2&\ge&C2^{-\delta n_0d(2\gamma -1)}\sum_zW_{k,z}(t)
^{2\gamma}\\
&\ge&C2^{-\delta n_0d(2\gamma-1)}\left[\sum_z\mathbf{1}(W_{k,z}(t)>0)\right]^
{1-2\gamma}V_{n_0}(t)^{2\gamma}.
\end{eqnarray*}
Therefore on $\mathcal{A}_{3,k}(\ell-1)$ and for 
$t\in[k\ell+\ell-1,(k+1)\ell]$, we have (note that $3^{d\gamma^{\ell-1}}\le c$)
\begin{eqnarray*}
H_t&\ge & C2^{-\delta n_0d(\gamma-1/2)}2^{n_0d\gamma^{\ell-1}(1/2-\gamma)}
V_{n_0}(t)^\gamma \\
&\equiv &AV_{n_0}(t)^\gamma.
\end{eqnarray*}
Thus we may apply Lemma \ref{a3} along with our choices of $\ell$ and 
$\delta$ (recall (\ref{ell2}) and (\ref{delta})) to conclude that
\begin{eqnarray*}
\lefteqn{P\left(\{V_{n_0}(t_{n_0+1})>0\}\cap\mathcal{A}_0(K)\cap
\left\{\sum_{k=0}^{\bar
m-1}\mathbf{1}(\mathcal{A}_{3,k}(\ell-1))>\frac{\bar m}{2}\right\}\right)}\\
&\le& P\bigg(\{V_{n_0}(t_{n_0+1})>0\}\cap\mathcal{A}_0(K)  \\
&&\cap\left\{\int_0^
{t_{n_0+1}}\mathbf{1}(H_t\ge AV_{n_0}(t)^\gamma)dt\ge\frac{\bar m}{2}
=\frac{t_{n_0+1}}{2\ell}\right\}\bigg)\\
&\le& CA^{-2}K^{2-2\gamma}2^{-n_0/2}\\
&=& C2^{n_0d(2\gamma -1)(\delta+\gamma^{\ell-1})-n_0/2}K^{2-2\gamma}\\
&\le & C 2^{n_0(d\frac{3}{2}\gamma^{\ell-1}-\frac{1}{2})}K^{2-2\gamma}\\
&\le & C2^{-n_0/4}K^{2-2\gamma}<\varepsilon/8
\end{eqnarray*}
for $n_0\ge n(\varepsilon, K)$.  This together with Lemma \ref{decrease}
completes the proof of Lemma \ref{n0}. \qed
\newline
We now turn to the
\newline
\noindent
\textbf{Proof of Lemma \ref{decrease}.}  Fix 
$0\le k< \bar m$.  For $0\le j<\ell$ let $\eta_{j,k}
=\sum_z\mathbf{1}(W_{k,z}(k\ell+j)>0)$ and
\[
\mathcal{A}_{0,k}(j,K)=\left\{\sup_{t\le k\ell+j}\sum_x u(t,x)\le K\right\}
\supset \mathcal{A}_0(K).
\]
By (\ref{lbH}) we may use Lemma \ref{a3} to see that for $1\le j <\ell$,
\begin{eqnarray*}
\lefteqn{E\left(\eta_{j,k}\mathbf{1}(\mathcal{A}_{0,k}(j,K))\bigg|\mathcal{F}
_{k\ell+j-1}\right)}\\
&\le&\sum_z\mathbf{1}(\mathcal{A}_{0,k}(j-1,K))P\left(W_{k,z}(k\ell+j)>0\bigg|
\mathcal{F}_{k\ell+j-1}\right)\\
&\le&\mathbf{1}(\mathcal{A}_{0,k}(j-1,K))\sum_zC2^{\delta n_0d(2\gamma-1)}
W_{k,z}(k\ell+j-1)^{2-2\gamma}.
\end{eqnarray*}
Now note that on $\mathcal{A}_{0,k}(j-1,K)$, we have $\sum_z W_{k,z}(k\ell+j-1)
\le K$ (by (\ref{wsum})).  We apply Lemma \ref{even} with 
\[
g(z)=W_{k,z}(k\ell+j-1) \hbox{ and }f(z)=\mathbf{1}(W_{k.z}(k\ell+j-1)>0)
\]
to see that on $\mathcal{A}_{3,k}
(j-1)$,
\begin{eqnarray*}
E\left(\eta_{j,k}\mathbf{1}(\mathcal{A}_{0,k}(j,K))\bigg|\mathcal{F}_{k\ell+j-1}
\right)&\le &C2^{2\delta n_0d(2\gamma-1)}K^{2-2\gamma}\eta_{j-1,k}^{2\gamma-1}
\\
&\le &C2^{n_0d(\gamma^{j-1}+\delta)(2\gamma-1)}K^{2-2\gamma}.\\
\end{eqnarray*}
Markov's inequality implies for $1\le j<\ell$,
\begin{eqnarray*}
\lefteqn{P\left(\mathcal{A}_{3,k}(j)^c\cap\mathcal{A}_{0,k}(j,K)\cap
\mathcal{A}_{3,k}(j-1)\bigg|\mathcal{F}_{k\ell+j-1}\right)}\\
&\le &\frac{C2^{n_0d(\gamma^{j-1}+\delta)(2\gamma-1)}K^{2-2\gamma}}
{3^{d\gamma^ j}2^{n_0d\gamma ^j}}\\
&\le &CK^{2-2\gamma}2^{n_0d(\delta(2\gamma-1)-\gamma^{j-1}(1-\gamma))}\\
&\le &CK^{2-2\gamma}2^{-n_0d\gamma^\ell(1-\gamma)/2}\equiv p_{n_0}',
\end{eqnarray*}
the last by (\ref{delta}).  Since there are at most $3^d2^{n_0d}$ sites $z$ 
in $D_{N_{n_0}}$, $P\left(\mathcal{A}_{3,k}(0)\right)=1$ and so from the above
we have, 
\begin{eqnarray*}
\lefteqn{P\left(\mathcal{A}_{3,k}(\ell-1)^c\cap\mathcal{A}_{0,k}(\ell-1,K)
      \bigg|\mathcal{F}_{k\ell}\right)}\\
&\le & P\left(\bigcup_{j=1}^{\ell-1}\left[\mathcal{A}_{3,k}(j)^c
\cap\mathcal{A}_{0,k}(j,K)\right]\bigg|\mathcal{F}_{k\ell}\right)\\
&\le &\sum_{j=1}^{\ell-1}P\left(\mathcal{A}_{3,k}(j)^c\cap\mathcal{A}_{0,k}
(j,K)\cap(\mathcal{A}_{3,k}(j-1)\cup\mathcal{A}_{0,k}(j-1,K)^c)\bigg|
\mathcal{F}_{k\ell}\right)\\
&=& \sum_{j=1}^{\ell-1}P\left(\mathcal{A}_{3,k}(j)^c\cap\mathcal{A}_{0,k}
(j,K)\cap\mathcal{A}_{3,k}(j-1)\bigg|\mathcal{F}_{k\ell}\right)\\
&\le &\ell p_{n_0}'\equiv p_{n_0}<\frac{1}{4},\\
\end{eqnarray*}
if $n_0>n(\ell,K,\gamma)$.  
Allowing $k<\bar m$ to vary, let
\[
d_k=\mathbf{1}(\mathcal{A}_{3,k}(\ell-1)^c\cap\mathcal{A}_{0,k}(\ell-1,K))\in
\mathcal{F}_{(k+1)\ell}, 
\]
and $M_n=\sum_{k=1}^{n-1}d_k-E(d_k|\mathcal{F}_{k\ell})$, $n<\bar m$.  Clearly
$(M_n,\mathcal{F}_{n\ell})_{n<\bar m}$ is a martingale and so if $n_0\ge
n(\ell,K,\gamma,\varepsilon)$, 
\begin{eqnarray*}
P(\sum_{k=0}^{\bar m -1}d_k\ge \bar m/2)&\le & P\left(M_{\bar m}\ge \bar m
\left({1\over2}-p_{n_0}\right)\right)\\
&\le & E\left(M_{\bar m}^2\right){\bar m}^{-2}\left({1\over 2}-p_{n_0}\right)^
{-2}\\
&\le &16{\bar m}^{-2}E\left(\sum_{k=1}^{\bar m -1}d_k^2\right)\\
&\le &16{\bar m}^{-1}p_{n_0}<{\varepsilon\over 8}.
\end{eqnarray*}
The probability we have to bound is no bigger than that on the left hand side
of the above equation and so the proof is complete. \qed

\ \newline
Now we can complete the 
\newline
\noindent
\textbf{Proof of Theorem \ref{t1}.}  Recall that 
\[
\mathcal{K}(t_\infty)\supset\mathcal{A}_0(K)
	\cap\left[\bigcap_{n=n_0}^\infty \mathcal{A}_{1,n}\right]
	\cap\left[\bigcap_{n=n_0}^\infty \mathcal{A}_{2,n}\right].
\]
Using (\ref{Kbound}), choose $K$ so large that 
\[
P\left(\mathcal{A}_0(K)^c\right)\le \frac{\varepsilon}{8}
\]
and then $n_0$ large enough so that all of the above bounds are valid.
Take complements in the above inclusion and consider the 
first value of $n$ so that $\omega\in \mathcal{A}_{1,n}^c$ or $\mathcal{A}_{2,
n}^c$ to see that (recall $\mathcal{A}_{i,n_0-1}$ is the entire space)
\begin{eqnarray*}
\mathcal{K}(t_\infty)^c\subset \mathcal{A}_0(K)^c&\cup&\left[\bigcup_{n=n_0}
^{\infty}\left(\mathcal{A}_{2,n}^c\cap\mathcal{A}_{1,n-1}\cap\mathcal{A}_0(K)
\right)\right]\\
&\cup&\left[\bigcup_{n=n_0}^{\infty}\left(\mathcal{A}_{1,n}^c\cap\mathcal{A}
_{2,n-1}\cap\mathcal{A}_0(K)\right)\right].
\end{eqnarray*}
Therefore Lemmas \ref{l.est-a}, \ref{a2} and \ref{n0} and our choice of $K$ 
imply
\begin{eqnarray*}
P(\mathcal{K}(t_\infty)^c)&\le&{\varepsilon\over 8}+{\varepsilon\over 4}
+\sum_{n=n_0+1}^\infty2^{n_0-n}{\varepsilon\over 4}+{\varepsilon\over 4}\\
&<&\varepsilon.
\end{eqnarray*}
This proves (\ref{main-est}), and finishes the proof of the theorem.
\qed

\section{Proof of Theorems \ref{3.4} and \ref{3.7}}
  \setcounter{equation}{0}

We first introduce a setting for mutually catalytic branching 
models.  Let $Q=(q_{xy})$ be the $Q$-matrix for a continuous time
$\mathbf{Z}^d$-valued Markov chain $\xi_t$ with semigroup $P_t$ and 
transition functions $\{p_t(x,y):t\ge 0, x,y\in\mathbf{Z}^d\}$.  If 
$|x|=|(x_1,\dots,x_d)|=\sum_{i=1}^d|x_i|$, $(x\in\mathbf{Z}^d)$, we assume the 
following hypotheses introduced in \cite{dp98}:
\begin{eqnarray*}
&(H1)& \Vert q\Vert_\infty=\sup_x|q_{xx}|<\infty.\\
&(H2)& \hbox{For each }x,y\in\mathbf{Z}^d,\ q_{xy}=q_{yx} \hbox{ and so }
p_t(x,y)=p_t(y,x).\\
&(H3)& \hbox{There are increasing positive functions }c(T,\lambda)\hbox{ and }
\lambda'(\lambda) \hbox{ such that}\\
&\ & \sum_y\left(|q_{xy}| + p_t(x,y)\right)\exp\left(\lambda|y|\right)\le
c(T,\lambda)\exp\left(\lambda'(\lambda)|x|\right)   \\
&&\qquad\forall t\in[0,T], x\in \mathbf{Z}^d.
\end{eqnarray*}     
These conditions are satisfied by a continuous time symmetric random
walk with subexponential tail (Lemma 2.1 of \cite{dp98}) and in particular
by the nearest neighbor random walk considered in the introduction for 
which $q_{xy}=\mathbf{1}(|x-y|=1)-2d\mathbf{1}(x=y)$.  Our generalized 
mutually catalytic system is then
\begin{eqnarray}
\label{system3}
U_t(x)&=&U_0(x)+\int_0^tQU_s(x)ds+\int_0^t\sqrt{U_s(x)V_s(x)} dB_{1,x}(s), 
 \nonumber \\
V_t(x)&=&V_0(x)+\int_0^tQV_s(x)ds+\int_0^t\sqrt{U_s(x)V_s(x)} dB_{2,x}(s), 
 	\\
x\in \mathbf{Z}^d,&\quad&  U_0, V_0\in M_F(\mathbf{Z}^d).\nonumber	
\end{eqnarray}
Here, $\{B_{i,x}(t)\}_{x\in\mathbf{Z}^d; i=1,2}$ is a collection of independent
$\mathcal{F}_t$-Brownian motions on some filtered probability space.
The weak existence and uniqueness of solutions in $C([0,\infty),M_F(\mathbf{Z}^d)^2)$
described in the introduction for $Q=\Delta$ 
continues to hold and we let $P_{U_0,V_0}$ continue to denote the 
unique law of the solution on this space of paths.

Theorem \ref{3.4} continues to hold without change in this more general
setting as we now show.\newline

\noindent 
\textbf{Proof of Theorem \ref{3.4}}. 
By Theorem 2.2(b)(ii) of \cite{dp98} $V_t(x)$ has mean $U_0P_t(x)$ and
variance
\begin{eqnarray*}
\lefteqn{\sum_y\int_0^t p_{t-s}(y,x)^2 U_0P_s(y) V_0P_s(y)ds}\\
&\leq& \int_0^t \left(\sum_{y_1}
p_{t-s}(y_1,x)U_0P_s(y_1)\right)\left(\sum_{y_2}p_{t-s}(y_2,x)V_0P_s(y_2)
\right)ds\\
&=&\int_0^t U_0 P_t(x)V_0P_t(x)ds =tU_0P_t(x)V_0P_t(x).
\end{eqnarray*}
By Chebychev's inequality
\begin{eqnarray*}
P_{U_0,V_0}\left( U_t(x)\leq \frac{1}{2}U_0P_t(x)\right) 
&\leq& P_{U_0,V_0}\left( |U_t(x)-U_0P_t(x)|\geq \frac{1}{2} 
	U_0P_t(x)\right)	\\
&\leq& 4t\frac{V_0P_t(x)}{U_0P_t(x)}.
\end{eqnarray*}
By (\ref{3.16}) if $\varepsilon >0$ and $t>t_0$, we may choose $x_0$ so that
$4tV_0P_t(x_0)/U_0P_t(x_0)<\varepsilon$ and so 
\[
P_{U_0,V_0}(\langle U_t,1\rangle >0)\geq P_{U_0,V_0}\left(
U_t(x_0)>\frac{1}{2} U_0P_t(x_0)\right) >1-\varepsilon .
\]
Since $\langle U_t,1\rangle$ is a non-negative martingale this shows that
\[
P_{U_0,V_0}\left(\langle U_t,1\rangle >0\quad\forall t >0\right) =1.
\]
The result follows by symmetry.\qed

It is easy to choose initial conditions satisfying (\ref{3.16}) for simple
symmetric random walk in $\mathbf{Z}^d$.

\begin{prop}
\label{3.5}
Assume $Q=\Delta$ so that $\{\xi_t\}$ is simple symmetric random walk on 
$\mathbf{Z}^d$ with jump rate $2d$.  Suppose there are $m>n$ in $\mathbf{Z}$ 
such that
\begin{eqnarray}
\label{3.17}
U_0\left( [m,\infty )\times\mathbf{Z} ^{d-1}\right) &=&0,\quad V_0\left(
[m,\infty )\times\mathbf{Z} ^{d-1}\right) >0\\
U_0\left( (-\infty ,n]\times\mathbf{Z} ^{d-1}\right)&>&0,\quad V_0\left(
(-\infty ,n]\times \mathbf{Z} ^{d-1}\right) =0.	\nonumber
\end{eqnarray}
Then (\ref{3.16}) holds and hence
\[
P_{U_0,V_0}\left(\langle U_t,1\rangle\langle V_t,1\rangle >0\quad\forall
t>0\right) =1.
\]
\end{prop}
We need an elementary estimate for simple random walk.

\begin{lemma}
\label{3.6} 
Let $\{\xi_t \}$ be simple symmetric random walk on
$\mathbf{Z}^d$. Then
\[
\frac{e^{-td}}{\prod\limits_{i=1}^d |x_i|!} t^{|x|}
\leq P_0(\xi_t =x)\leq \frac{1}{\prod_{i=1}^d |x_i|!} t^{|x|}
\]
where we recall that $|x|=\sum\limits_1^d|x_i|$.
\end{lemma}

\noindent \textbf{Proof}. Suppose $d=1$ and $\xi_t$ jumps with rate $\lambda >0$.
Then for $x\geq 0$,
\begin{eqnarray}
\label{3.18}
p_t^{(\lambda )}(x) &\equiv& P_0(\xi_t =x)	\\
&=& \sum_{n=0}^\infty
P(\xi\hbox{ has $n+x$ steps to the right up to time $t$}) \nonumber\\
&=& \sum_{n=0}^\infty {2n+x\choose n+x} 2^{-2n-x}e^{-\lambda t}
\frac{(\lambda t)^{2n+x}}{(2n+x)!} \nonumber\\
&=& e^{-\lambda t} (\lambda t/2)^x \sum\limits_{n=0}^\infty 
\frac{(\lambda t/2)^{2n}}{n!(n+x)!} \nonumber\\
&\equiv& e^{-\lambda t} (\lambda t/2)^x \sigma (x).	\nonumber
\end{eqnarray}
Clearly
\begin{eqnarray*}
\frac{1}{x!}\leq \sigma (x)&\leq& \frac{1}{x!}
\sum_{n=0}^\infty {2n\choose n}2^{-2n}\frac{(\lambda t)^{2n}}{(2n)!}\\
&\leq& \frac{1}{x!}\sum_{n=0}^\infty \frac{(\lambda t)^{2n}}{(2n)!}\leq
\frac{e^{\lambda t}}{x!}.
\end{eqnarray*}
Put this into (\ref{3.18}) and use symmetry in $x$ to get
\begin{equation}
\label{3.19}
e^{-\lambda t} \frac{(\lambda t/2)^{|x|}}{|x|!}\leq p_t^{(\lambda)}(x)\leq
\frac{(\lambda t/2)^{|x|}}{|x|!}\qquad\forall x\in\mathbf{Z} .
\end{equation}
Since $P_0(\xi_t =x)=\prod_{i=1}^d p_t^{(1)}(x_i)$, the result
follows.\qed \newline

\noindent \textbf{Proof of Proposition \ref{3.5}}. Choose $v\in \mathbf{Z}^d$ 
with $v_1\geq m$
such that $V_0(v)>0$. Then (\ref{3.17}) shows that for $t>2$, 
$x=(x_1,v_2,\ldots ,v_d)$ and $x_1>v_1$,
\begin{eqnarray*}
\frac{U_0P_t(x)}{V_0P_t(x)} 
&\leq& \sum_{y} \mathbf{1}(y_1<m)U_0(y)\frac{t^{\sum\limits_1^d
(|y_i-x_i|-|v_i-x_i|)}}{\prod\limits_1^d |y_i-x_i|! V_0(v)e^{-dt}}
\prod\limits_1^d |v_i-x_i|!\\
&\leq& \frac{e^{2dt}}{V_0(v)}(x_1-v_1)! \sum_{y_1<m}
\frac{t^{v_1-y_1}}{|y_1-x_1|!} \sum_{y_2,\ldots ,y_d} U_0(y)
\frac{t^{\sum_2^d |y_i-x_i|}}{\prod_2^d |y_i-x_i|!}\\
&\leq& \frac{e^{3dt}}{V_0(v)}\sum_{y_1<m}
\frac{t^{(v_1-m)+(m-y_1)}}{(x_1-y_1)\ldots (x_1-v_1+1)} U_0^{(1)}(y_1),\\
\end{eqnarray*}
where $U_0^{(1)}(y_1)=\sum\limits_{y_2\ldots y_d}U_0(y_1,y_2,\ldots ,y_d)$
is the first marginal of $U_0$. Setting $k=m-y_1$, we see that
\begin{eqnarray*}
\frac{U_0P_t(x)}{V_0P_t(x)} 
&\leq& t^{v_1-m} \frac{e^{3dt}}{
V_0(v)}\sum_{k=1}^\infty \frac{t^k}{k!}\cdot\frac{k!}{(x_1-m+k)\ldots
(x_1-v_1+1)} \\
&&\cdot \, U_0^{(1)}(m-k)\\
&\leq& c(t,v,m)e^{(3d+1)t}\sum_{k=1}^\infty \frac{U_0^{(1)}(m-k)}{
(x_1-v_1+1)}\\
&\leq& \frac{c'(t,v,m)\langle U_0,1\rangle}{x_1-v_1+1}\\
&\rightarrow& 0\quad\hbox{as}\quad x_1\rightarrow \infty .
\end{eqnarray*}
This, and a symmetrical argument for the reciprocal, establish (\ref{3.16})
and we are done.\qed

Next, we turn to the main task of this section, proving Theorem
\ref{3.7} in our more general Markov chain setting.  In this case 
we need to add a pair of hypotheses.
\begin{theorem}
\label{3.7q} Assume 
\begin{equation}
\label{4.3}
(H3)\hbox{ holds with }\lambda'(\lambda)=\lambda
\end{equation}
and
\begin{equation}
\label{4.4}
\inf_{s\le t,x\in \mathbf{Z}^d}p_s(x,x)=\varepsilon_0(t)>0
\end{equation}
Under the hypotheses on $V_0$ in Theorem \ref{3.7}, the conclusion of that
result holds.
\end{theorem}
\textbf{Remark.} The hypotheses added above hold for any continuous time 
random walk with subexponential jump distributions, as (\ref{4.4}) is 
trivial and (\ref{4.3}) is proved in Lemma 2.1 of \cite{dp98}.

\noindent
\textbf{Notation.} $\phi_\lambda(x)=e^{\lambda|x|}$ for
$x\in\mathbf{Z}^d$.
\begin{lemma}
\label{3.8} Assume the hypotheses of Theorem \ref{3.7q}. Then
$\forall\lambda >0$, $\varepsilon >0$, $T>0$ there is a $C_{T,\varepsilon
,\lambda}>0$ such that
\[
\varepsilon_0 (T)\phi_{-\lambda}(x)\leq P_t\phi_{-\lambda}(x)\leq
C_{T,\varepsilon ,\lambda}\phi_{-\lambda +\varepsilon} (x)\qquad\forall 
t\leq T.
\]
\end{lemma}

\noindent \textbf{Proof}.  For the lower bound, observe that
$P_t\phi_{-\lambda}(x)\geq p_t(x,x)\phi_{-\lambda}(x)\geq
\varepsilon_0 (T)\phi_{-\lambda}(x)$ for all $t\leq T$ by (\ref{4.4}).

For the upper bound, note that $(H_2)$ implies
\[
p_t(x,y)=p_t(y,x)\leq C_{T,\lambda -\varepsilon }e^{(\lambda -\varepsilon
)|y|-(\lambda -\varepsilon )|x|}\quad\hbox{for}\quad t\leq T
\]
and so for $t\leq T$
\[
P_t(\phi_{-\lambda} )(x)\leq C_{T,\lambda -\varepsilon}\sum_{y\in \mathbf{Z}
^d} e^{-\varepsilon |y|}e^{-(\lambda -\varepsilon )|x|}\leq C_{T,\varepsilon
,\lambda} \phi_{-(\lambda -\varepsilon )}(x).
\]
This finishes the proof of Lemma \ref{3.8}.
\qed \newline

\noindent \textbf{Proof of Theorem \ref{3.7q}}. The hypotheses on
$\{\lambda_i\}$ allow us to choose
$\beta$ such that 
\begin{equation}
\label{3.21}
2\lambda_1 -\Big(\frac{\lambda_0+\lambda_2}{2}\Big)<\beta 
<\lambda_2
\end{equation}
and then $\alpha$ such that
\begin{equation}
\label{3.22}
2\lambda_1 -\beta <\alpha <(\lambda_0+\lambda_1 )/2.
\end{equation}
For now we fix $\eta\in (0,1]$ (recall $U_0\in\eta\phi_{-\lambda_0}$)and 
will specify its value later in the proof. It is easy to modify the 
derivation of Theorem 2.2(c) of \cite{dp98}  to see that 
\begin{equation}
\label{3.23}
\langle V_t,\phi_\beta\rangle =\langle V_0,\phi_\beta \rangle
+\int_0^t \langle V_s,Q\phi_\beta\rangle  ds+M_t^V (\phi_\beta),
\end{equation}
where $M_t^V(\phi_\beta )$ is a continuous square integrable martingale
such that 
\[
\langle M^V(\phi_\beta )\rangle_t=\int_0^t \langle
U_sV_s,\phi_\beta ^2\rangle ds.
\]
To see this note that $\beta <\lambda_2$ shows that $\langle V_0,
\phi_\beta\rangle <\infty$, and that (\ref{4.3}) implies 
\begin{equation}
\label{3.24}
|Q\phi_\beta |\leq c(\beta )\phi_\beta
\end{equation}
so that
\begin{eqnarray}
\label{3.25}
E\left( \langle V_s,|Q\phi_\beta |\rangle \right)
&\leq& c(\beta ) \langle V_0,P_s\phi _\beta \rangle \\
&\leq& c(\beta ) C_{T,\beta } \langle V_0, \phi_\beta \rangle\qquad
\hbox{(by (\ref{4.3}) )} \nonumber\\
&<& \infty .\nonumber
\end{eqnarray}
In addition we use the fact that
\begin{eqnarray*}
E\left( \int_0^t \langle U_sV_s, \phi_\beta^2 \rangle  ds\right)
&=& \int_0^t \sum_x e^{2\beta |x|} U_0P_s(x) V_0 P_s(x) ds\\
&\ &\qquad\qquad\hbox{(Theorem 2.2 of \cite {dp98})}\\
&\leq& \int_0^t \sum_x e^{2\beta |x|} c_{t,\varepsilon
,\lambda_0} \eta e^{-(\lambda_0-\varepsilon )|x|}c_{t,\varepsilon
,\lambda_2}c_2e^{-(\lambda_2-\varepsilon )|x|} ds\\
&\ &\qquad\qquad\hbox{(Lemma \ref{3.8})}\\
&<& \infty\qquad\hbox{for any}\quad t>0,
\end{eqnarray*}
for $\varepsilon >0$ small enough because $2\beta <2\lambda_2 <\lambda_0
+\lambda_2$. (\ref{3.25}) and the above show that
\begin{equation}
\label{3.26}
V_T^{(\beta )} \equiv \sup_{t\leq T}\langle V_t,\phi_\beta
\rangle\in L^1 \hbox{ and  }\sup_{0<\eta \leq
1}P_{U_0,V_0}\left( V_T^{(\beta )}\right) \leq C_T^V <\infty .
\end{equation}
A similar argument (now use $\alpha <(\lambda_0 +\lambda_1)/2<\lambda_0$)
shows that
\begin{equation}
\label{3.27}
U_T^{(\alpha )}=\sup_{t\leq T}\langle U_t,\phi_\alpha\rangle\in \mathbf{L}^1
\hbox{ and  }\sup_{0<\eta\leq 1}
P_{U_0,V_0}\left( U_T^{(\alpha )}\right) \leq C_T^U<\infty .
\end{equation}
Now fix $\varepsilon >0$. We claim there is a $t_0>0$, independent of the
choice of $\eta \in (0,1]$, such that
\begin{equation}
\label{3.28}
P_{U_0,V_0} \left( V(t,x)\geq \frac{1}{2} V_0P_t(x)\quad\forall t\leq
t_0\quad \forall x\in \mathbf{Z}^d\right) >1-\varepsilon .
\end{equation}
If $N(t,x)=V(t,x)-V_0P_t(x)$, then Theorem 2.2(b) of \cite{dp98} shows that if
$0\leq t\leq u\leq 1$ then $N(u,x)-N(t,x)=N_{t,x}^{(1)}
(u)+N_{u,x}^{(2)}(t)$, where
\begin{eqnarray*}
N_{t,x}^{(1)} (u)&=& \sum_y \int_t^u
p_{u-s}(y,x)\sqrt{U_s(y)V_s(y)} dB_{2,y}(s)\\
\noalign{\hbox{and}}
N_{u,x}^{(2)}(t) &=& \sum_y \int_0^t \big(
p_{u-s}(y,x)-p_{t-s}(y,x)\big)\sqrt{U_s(y)V_s(y)} dB_{2,y}(s).
\end{eqnarray*}
Our hypotheses on $\{\xi _t\}$ imply that
\begin{equation}
\label{3.29}
|p_r(x,y)-p_s(x,y)|=\left| \int_s^r Qp_w(\cdot ,y)(x) dw\right|
\leq \| q\|_\infty |r-s|.
\end{equation}
We have for any $\delta >0$,
\begin{eqnarray}
\label{3.30}
\langle N_{u,x}^{(2)}\rangle_t 
&=& \sum_y \int_0^t \left( p_{u-s}(y,x)-p_{t-s}(y,x)\right)
^2 U_s(y)V_s(y) ds\\
&\leq& \| q\|_\infty (u-t)\int_0^t \sum_y \big(
p_{u-s}(y,x)+p_{t-s}(y,x)\big)   \nonumber\\
&& \cdot\, U_1^{(\alpha )} V_1^{(\beta)}
\phi_{-\alpha}(y)\phi_{-\beta }(y) ds	\nonumber\\
&&\hskip1.6in \hbox{(by (\ref{3.29}), (\ref{3.26}) and (\ref{3.27}))}
		\nonumber\\
&\leq& U_1^{(\alpha)} V_1^{(\beta)} \| q\|_\infty (u-t)
C_{1,\delta,\alpha+\beta}
\phi_{-\alpha -\beta +\delta} (x)\quad\hbox{(Lemma \ref{3.8})}
\nonumber\\
&\equiv& c(\delta ) U_1^{(\alpha)} V_1^{(\beta)} (u-t) \phi_{-\alpha-\beta
+\delta } (x)\nonumber
\end{eqnarray}
and, using similar reasoning,
\begin{eqnarray*}
{\langle N_{t,x}^{(1)}\rangle }_u
&=& \sum_y \int_t^u p_{u-s}(y,x)^2 U_s(y)V_s(y) ds\\
&\leq& U_1^{(\alpha)} V_1^{(\beta)} \int_t^u \sum_y
p_{u-s}(x,y)\phi_{-\alpha} (y)\phi_{-\beta }(y) ds\\
&\leq& U_1^{(\alpha)} V_1^{(\beta)} (u-t)\sup_{s\leq u} P_s\phi_{-\alpha
-\beta}(x)\\
&\leq& c(\delta ) U_1^{(\alpha)} V_1^{(\beta)} (u-t)\phi_{-\alpha -\beta
+\delta}(x).
\end{eqnarray*}
Choose $\delta >0$ small enough so that (see (\ref{3.22}))
\begin{equation}
\label{3.31}
2\lambda \equiv \alpha+\beta -\delta >2\lambda _1.
\end{equation}
Let $\Delta (n,x)=(n+1)^{1/2}2^{-n/2}(|x|+1)^{1/2} \phi_{-\lambda}(x)$,
$x\in \mathbf{Z}^d$, $n\in \mathbf{N}$. Then for $K,K_1>0$
\begin{eqnarray*}
\lefteqn{P_{U_0,V_0}
\left( 
    \left| N \left(\frac{j+1}{2^n},x\right) -N
        \left( \frac{j}{2^n},x\right) 
    \right| 
\geq K\Delta (n,x),\quad U_1^{(\alpha)} V_1^{(\beta)} \leq K_1^2
\right)}\\
&\leq& P_{U_0,V_0} 
\bigg( 
    \left| N_{j/2^n,x}^{(1)} \left( (j+1)/2^n\right) 
    \right|
      \geq \frac{K}{2}\Delta (n,x), \\
&&\qquad\langle N_{j2^{-n},x}^{(1)}\rangle
       \left( (j+1)2^{-n}\right)\leq c(\delta )K_1^22^{-n}\phi_{-2\lambda}(x)
\bigg) \\
&& + P_{U_0,V_0}
\bigg(
    \left| N_{(j+1)2^{-n},x}^{(2)} 
       \left( j2^n\right) 
    \right|
       \geq \frac{K}{2}\Delta (n,x), \\
&&\qquad\langle N_{(j+1)2^{-n},x}^{(2)}\rangle
       (j2^{-n})\leq c(\delta )K_1^22^{-n}\phi_{-2\lambda}(x)
\bigg) \\
&\leq& 2P
\left(
    \sup_{s\leq 1} |B_s| >\frac{1}{2} K\Delta (n,x)
        {\big( c(\delta )K_1^2 2^{-n}\phi_{-2\lambda}(x)\big)}^{-1/2}
\right) ,
\end{eqnarray*}
where $B$ is a one-dimensional Brownian motion and we have used the
Dubins-Schwarz Theorem in the last line. An elementary estimate on the
Gaussian tail and the fact that
\begin{eqnarray*}
\lefteqn{\frac{1}{2} K\Delta (n,x){\left( c(\delta )K_1^22^{-n}
\phi_{-2\lambda}(x)\right)}^{-1/2}}  \\
&&= KK_1^{-1} (n+1)^{1/2}{\left(
|x|+1\right)}^{1/2} \left( 2\sqrt{c \left( \delta \right) }\right)^{-1}
\end{eqnarray*}
shows that if we set $L=K^2/4K_1^2 c(\delta )$ and assume $L\geq 1$, then we
have
\begin{eqnarray*}
\lefteqn{P_{U_0,V_0}
\bigg(
   \left| N\left(\frac{j+1}{2^n},x\right) -N\left(\frac{j}{2^n},x\right)
   \right| \geq K\Delta (n,x)}	\\
&&\qquad\hbox{for some}\quad 0\leq j<2^n, x\in
\mathbf{Z}^d\quad\hbox{and}\quad n\in \mathbf{N} 
\bigg) \\
&\leq& 8\sum_{n=1}^\infty \sum_{x\in \mathbf{Z}^d}2^n\exp
\left( \frac{-K^2(n+1)(|x|+1)}{4K_1^2c(\delta )}\right)
+ P_{U_0,V_0} \left( U_1^{(\alpha)} V_1^{(\beta)} >K_1^2\right)\\
&\leq& 4\sum_{x\in \mathbf{Z}^d} \sum_{n=1}^\infty \exp
   \left\{ -\left[ \frac{K^2(|x|+1)}{4K_1^2c(\delta )}-\log 2\right]
    (n+1)\right\}	\\
&& + P_{U_0,V_0} \left( U_1^{(\alpha)} >K_1\right) +P_{U_0V_0}\left( V_1^{(\beta)}
            >K_1\right) \\
&\leq& c\sum_{x\in \mathbf{Z}^d} \exp \left\{ -L(|x|+1)\right\}
+(K_1)^{-1}\left( C_1^U+C_1^V\right),
\end{eqnarray*}
where we have used (\ref{3.26}) 
and (\ref{3.27}) in the last line.
First choose $K_1$ and then $K$ sufficiently large so that the above
expression is less than $\varepsilon$ (and $L\ge 1$). Note that the choice 
of $K_1$ and $K$ may be made independently of $\eta\in (0,1]$. Therefore 
off a set of $P_{U_0,V_0}$-measure at most $\varepsilon$ if 
$2^{-n_0}\leq t< 2^{1-n_0}$ $(n_0\in \mathbf{N})$ and 
$t=\sum_{n_0}^\infty j_n2^{-n}$ where $j_{n_0}=1$ and $j_n\in \{ 0,1\}$ for 
$n>n_0$, then for all $x$ in $\mathbf{Z}^d$
\begin{eqnarray*}
|N (t,x)| &\leq& \left(\sum_{n=n_0}^\infty
j_n2^{-n/2}(n+1)^{1/2}\right) K{(|x|+1)}^{1/2} \phi_{-\lambda }(x)\\
&\leq& cK{(t\log 1/t)}^{1/2}{(|x|+1)}^{1/2}\phi_{-\lambda}(x)\varepsilon_0
(1)^{-1} c_1^{-1}\phi_{\lambda_1}(x) P_tV_0(x)\\
&&\hskip1in\hbox{(by the lower bound in Lemma \ref{3.8})}\\
&\leq& cK{(t\log 1/t)}^{1/2} P_tV_0(x)\qquad\hbox{(since $\lambda_1
<\lambda$).}
\end{eqnarray*}
Hence we may choose $t_0>0$ sufficiently small (independent of $\eta\in
(0,1]$ and $c_0$) such that
\[
P_{U_0,V_0}\left( |N(t,x)|\leq \frac{1}{2} P_tV_0(x)\qquad\hbox{for
all}\quad 0\leq t\leq t_0, x\in \mathbf{Z}^d\right) >1-\varepsilon .
\]
This proves (\ref{3.28}).

Let 
\[
\sigma = \inf \left\{ t: V_t(x)<\frac{1}{2}P_tV_0(x)\quad\hbox{for
some}\quad x\in \mathbf{Z}^d \right\} \wedge 1.
\]
As for (\ref{3.23}) we have
\begin{equation}
\label{3.32}
\langle U_t,\phi_\alpha \rangle = \langle U_0,\phi_\alpha\rangle
+\int_0^t \langle U_s,Q\phi_\alpha \rangle  ds+M_t^U(\phi_\alpha)
\end{equation}
where $M_t^U(\phi_\alpha )$is a continuous square-integrable martingale
such that for $t\le \sigma$,
\begin{eqnarray*}
\frac{d}{dt}{\langle M^U(\phi_\alpha )\rangle} 
&=& {\langle U_tV_t,\phi_\alpha ^2\rangle}\\
&\geq& \frac{1}{2}\langle U_tP_tV_0,\phi_\alpha
^2\rangle\\
&\geq& \frac{c_1}{2} \varepsilon_0(1) \langle U_t,\phi_\alpha\rangle.
\end{eqnarray*}
In the last line we use Lemma \ref{3.8} and 
the fact that $\alpha >\lambda_1$ by (\ref{3.21}) and (\ref{3.22}).
Therefore
\begin{equation}
\label{3.33}
\frac{d}{dt} {\langle M^U(\phi_\alpha )\rangle}_t \geq c_{\ref{3.33}}\langle
U_t,\phi_\alpha \rangle\qquad\hbox{for}\quad t\leq\sigma .
\end{equation}
Let
\begin{eqnarray*}
C(t)&=&\int_0^t \langle U_s,\phi_\alpha \rangle ^{-1}d{\langle
M^U(\phi_\alpha )\rangle }_s,\\
\tau (t) &=& \inf\left\{ u:C(u)>t\right\}\qquad (\inf\phi =\infty ),\\
\noalign{\hbox{and}}
\tilde U (t) &=& \langle U_{\tau
(t)},\phi_\alpha\rangle\qquad\hbox{for}\quad t<C(\infty ).
\end{eqnarray*}
Then for $t<C(\sigma )$
\begin{eqnarray*}
\tilde U (t)&=& \langle U_0,\phi_\alpha\rangle +\int_0^{\tau
(t)}\langle U_s,Q\phi_\alpha\rangle ds+\tilde M (t)\\
&=& \langle U_0,\phi_\alpha\rangle +\int_0^t\langle
U_{\tau_r},Q\phi_\alpha\rangle \tau '(r) dr+\tilde M (t)
\end{eqnarray*}
where $\tilde M (t)=M_{\tau (t)\wedge\sigma}^U(\phi_\alpha )$ is a
continuous local martingale such that for $t<C(\sigma )$
\[
\frac{d}{dt} {\langle \tilde M \rangle }_t =\langle M^U(\phi_\alpha )\rangle
'\big(\tau (t)\big) \tau '(t)=\langle M^U(\phi_\alpha )\rangle '\big(\tau
(t)\big) C'\big(\tau (t)\big)^{-1}= \tilde U (t).
\]
Hence by enlarging the probability space if necessary we may assume there
is a filtration $(\tilde\mathcal{F}_t)$ and an $(\tilde\mathcal{F}_t)$-Brownian motion,
$B(t)$, such that $U_{\tau (t)}$, $\tilde U _t$ and 
$\tau '(t)\mathbf{1}(t<C(\sigma))$ are $(\tilde\mathcal{F}_t)$-adapted, 
$C(\sigma )$ is an $(\tilde\mathcal{F}_t)$-stopping time and 
\begin{equation}
\label{3.34}
\tilde U (t)=\langle U_0,\phi_\alpha\rangle +\int_0^t \langle
U_{\tau_r},Q\phi_\alpha\rangle \tau '(r) dr+\int_0^t\sqrt{\tilde U
(r)} dB_r,\ t<C(\sigma ).
\end{equation}
(\ref{3.33}) implies that $\tau '(r)\leq c_{\ref{3.33}}^{-1}$ 
for $r<C(\sigma )$ and
the analogue of (\ref{3.24}) now shows that for $r<C(\sigma )$,
\begin{eqnarray*}
\left|\langle\tilde U _r,Q\phi_\alpha\rangle\tau '(r)\right| 
&\leq& c(\alpha )c_{\ref{3.33}}^{-1}\tilde U _r\\
&\equiv& c_{\ref{3.35}} \tilde U _r.
\end{eqnarray*}
Let $\hat U_t$ be the pathwise unique solution of 
\begin{equation}
\label{3.35}
\hat U(t)=\langle U_0,\phi_\alpha\rangle +\int_0^t c_{\ref{3.35}} \hat U
_r dr+\int_0^t \sqrt{\hat U (r)} dB(r).
\end{equation}
A comparison theorem for stochastic differential equations (see Rogers and
Williams \cite{rw87}, V.43.1) shows that
\[
\tilde U (t)\leq \hat U (t)\qquad\hbox{for}\quad t<C(\sigma )
\]
and so
\begin{equation}
\label{3.36}
\langle U_t,\phi_\alpha\rangle\leq\hat U \big(
C(t)\big)\qquad\hbox{for}\quad t<\sigma .
\end{equation}
(\ref{3.33}) shows that $C(t)\geq c_{\ref{3.33}}t$ for 
$t\leq\sigma$ and so if $P_x$ is
the law of $\hat U$ starting at $x$ we have for $0<t_1\leq t_0$ (here
$\varepsilon >0$ is fixed and $t_0$ is as in (\ref{3.28}))
\begin{eqnarray*}
P_{U_0,V_0}(U_t=0\quad\forall t\geq t_1)
&=& P_{U_0,V_0} (U_{t_1}=0)\\
&\geq& P\left(
            \hat U \big( C(t_1)\big) =0,\quad t_1<\sigma
        \right) \qquad\hbox{(by (\ref{3.36}))}\\
&\geq& P_{\eta\langle\phi_{-\lambda_0},\phi_\alpha\rangle }
   \left(\hat U (c_{\ref{3.33}}t_1) =0\right) 
      -P_{U_0,V_0}(\sigma\leq t_0)\\
&\geq& P_1\left( \hat U (c_{\ref{3.33}}t_1)=0\right)
^{\eta\langle\phi_{-\lambda_0},\phi_\alpha\rangle } -\varepsilon \hbox{
(by (\ref{3.28}))},
\end{eqnarray*}
where we have used the multiplicative property of the superprocess $\hat U$
and the choice of $t_0$. Now $P_1 \big( \hat U (c_{\ref{3.33}} t_1) = 0 \big)
>0$ (in fact it is easy to get an explicit expression for this probability,
or alternatively one may use Girsanov's theorem and the fact that this
probability is positive if $c_{\ref{3.35}}=0$ in (\ref{3.35}).  Therefore
for $\eta>0$ sufficiently small the above 
probability is at least $1-2\varepsilon$. This completes the proof of 
Theorem \ref{3.7q}.\qed

\newcommand{\etalchar}[1]{$^{#1}$}
\providecommand{\bysame}{\leavevmode\hbox to3em{\hrulefill}\thinspace}


\end{document}